\documentclass[12pt]{amsart}

\usepackage[bookmarksopen,bookmarksdepth=2]{hyperref}
\usepackage{etex}
\usepackage{graphicx}
\usepackage{amssymb,amsmath,amsfonts,amsthm,epsfig,amscd,stmaryrd,mathrsfs }
\usepackage{csquotes}
\usepackage{stmaryrd}
\usepackage[all,cmtip,poly]{xy}
\usepackage{a4wide}
\usepackage{color}

\def\pp{\mathfrak{p}}

\def\AA{\mathcal{A}}
\def\SL{\mathrm{SL}}
\def\RGR{\mathscr{G}\kern-0.1em{\mathscr{R}}}
\def\Rfl{R^{\mathrm{fl}}}

\def\RR{\mathcal{R}}
\def\Rglob{R^{\mathrm{glob}}}
\def\Sym{\mathrm{Sym}}
 \def\Fbar{\overline{\F}}
\def\a{\mathfrak{a}}
 \def\Zbar{\overline{\Z}}
 \def\Rloc{R^{\mathrm{loc}}}
 \def\Gal{\mathrm{Gal}}
 \def\ad{\mathrm{ad}}
\def\T{\mathbf{T}}
\def\F{\mathbf{F}}
\def\Z{\mathbf{Z}}

\def\OL{\mathcal{O}}
\def\Frob{\mathrm{Frob}}

\def\Q{\mathbf{Q}}

\def\GL{\mathrm{GL}}
\def\Qbar{\overline{\Q}}
\def\rhobar{\overline{\rho}}
\def\m{\mathfrak{m}}
\def\A{\mathbf{A}}
\def\C{\mathbf{C}}

\numberwithin{equation}{section}

\begin{document}

\title{Reciprocity in the Langlands program since Fermat's Last Theorem}

\author{Frank Calegari}

\thanks{The  author was supported in part by NSF Grant
  DMS-2001097}

\address{The University of Chicago,
5734 S University Ave,
Chicago, IL 60637, USA  \texttt{fcale@uchicago.edu}}

\begin{abstract} 
The \emph{reciprocity conjecture} in the Langlands program links motives to automorphic forms.
The proof of Fermat's Last Theorem by Wiles~\cite{W,TW} introduced new tools to study reciprocity.
This survey reports on developments using these ideas (and their generalizations) in the last three decades.
\end{abstract}

\maketitle

\section{Introduction}

The \emph{reciprocity conjecture} in the Langlands program predicts a relationship
between pure motives\footnote{Here (in light of the standard conjectures~\cite{MR1265519}) one may take pure motives  up to numerical or homological equivalence.
 Conjecturally, one can also substitute
 (for irreducible motive) the notion of an irreducible weakly compatible system
of Galois representations~\cite{TaylorICM} or an irreducible geometric Galois representation in the sense of Fontaine--Mazur~\cite{FM}.}
 and automorphic representations.  The simplest version (as formulated  by Clozel~\cite[Conj~2.1]{Clozel}) states that
there should be a bijection
between irreducible motives~$M$ over a number field~$F$ with coefficients in~$\Qbar$ and
cuspidal algebraic representations~$\pi$ of~$\GL_n(\A_F)$ satisfying a number of explicit additional compatibilities,
including the equality of algebraic and analytic~$L$-functions~$L(M,s) = L(\pi,s)$. In light of multiplicity one theorems~\cite{MR623137}, this pins down
the correspondence uniquely.
There is also a version of this conjecture
for more general reductive groups, although its formulation requires some care (as was done by Buzzard and Gee~\cite{BG}).
Beyond the spectacular application by Wiles to Fermat's Last Theorem~\cite[Theorem~0.5]{W}, 
the Taylor--Wiles method~\cite{W,TW} 
 gave a completely new technique --- and to this date
the most successful one --- for studying the problem of reciprocity. The ideas in these two papers have sustained progress in the field
for almost\footnote{Wiles in~\cite{W} dates the completion of the proof to September~19, 1994.}
 30 years.
In  this  survey, we explain how the Taylor--Wiles method has evolved over this period
and where it stands today. \textbf{One warning:}   the intended audience for this document is 
entirely complementary to the audience for my talk ---
I shall assume more than a passing familiarity with the arguments of~\cite{W,TW}.
Moreover, this survey is as much a personal and historical\footnote{A whiggish history, naturally.
Even with this caveat, it should be clear  
that
the narrative arc of progress presented here at best represents my own interpretation of events. I have added a few quotes from first hand sources when I felt
they conveyed a sense of what the experts were thinking in a manner not easily obtainable from other sources.
For other survey articles on similar topics, see~\cite{BuzzardE,BreuilE}.}
discussion as a mathematical one ---  giving anything more than hints
on even a fraction of what is discussed here would be close to impossible given the space constraints and the competence of the author.
Even with the absence  of any real mathematical details in this paper, the sheer amount of activity in this field has led me to discard any discussion
of advances not directly related to~$R=\T$ theorems, which necessitates the omission of a lot of closely related beautiful  mathematics.

\subsection{The Fontaine--Mazur Conjecture}

Let~$F$ be a number field. The Fontaine--Mazur conjecture\footnote{Fontaine told me
(over a \emph{salad de g\'{e}siers} in Roscoff in 2009) that he and Mazur formulated their conjecture in the mid-80s. 
(Colmez pointed me towards these notes~\cite{FontaineNotes} from a talk given by Fontaine at the 1988 Mathematische Arbeitstagung in Bonn.)
He noted that Serre had originally been skeptical,
particularly of the claim  that any everywhere unramified representation inside~$\GL_n(\Qbar_p)$
must have finite image, and  set off to find a counterexample (using the construction of Golod--Shavarevich~\cite{GS}).
He (Serre) did not succeed!}~\cite{FM} predicts that any  continuous irreducible~$p$-adic Galois representation
\begin{equation*}
\rho: G_F \rightarrow \GL_n(\Qbar_p)
\end{equation*}
which is both unramified outside finitely many primes and potentially semi-stable (equivalently, de Rham~\cite{ColmezB}) at all places~$v|p$ should be associated to a motive~$M/F$ with coefficients in~$\Qbar$.
Any such~$\rho$ is automatically conjugate to a representation in~$\GL_n(E)$ for some finite field~$E/\Q_p$ and further stabilizes an~$\OL_E$-lattice. The 
 corresponding 
residual representation~$\rhobar: G_F \rightarrow \GL_n(k)$ where~$k = \OL_E/\pi_E$ is the residue field of~$E$ is unique up to semi-simplification.
 Let us assume here for expositional convenience that~$\rhobar$  is absolutely irreducible.
Following Mazur~\cite{Mazur}, one may define a universal deformation ring which parameterizes
all deformations of~$\rhobar$ unramified outside a finite set~$S$. 
One can then further impose local conditions to define deformation rings~$R$ whose~$\Qbar_p$-valued
points are associated to Galois representations which are de Rham at~$v|p$ with fixed Hodge--Tate weights.
Assuming the Fontaine--Mazur conjecture, these~$\Qbar_p$-valued points correspond to all pure motives~$M$ unramified outside~$S$ 
whose~$p$-adic realizations are  Galois representations with the same  local conditions at~$p$ and the same 
fixed residual representation~$\rhobar$.
Assuming the reciprocity conjecture, these motives should then be associated to a finite dimensional space of automorphic forms. This  leads to the extremely
non-trivial prediction that~$R$ has finitely many~$\Qbar_p$-valued points. 
The problem of reciprocity is now to link these~$\Qbar_p$-valued points of~$R$ to automorphic forms.

\subsection{\texorpdfstring{$R = \T$}{R eq T} theorems}

Associated to the (conjectural) space of automorphic forms corresponding  to~$\Qbar_p$-valued points of~$R$
is a ring of endomorphisms generated by Hecke operators. The na\"{\i}ve version of~$\T$ is defined to be the completion of this ring with respect to a maximal ideal~$\m$ defined in terms of~$\rhobar$. The mere existence of~$\m$ is itself conjectural, and amounts --- in the special case of odd absolutely irreducible~$2$-dimensional representations~$\rhobar$ of~$G_{\Q}$ --- to
 Serre's conjecture~\cite{MR885783}. Hence,  in the Taylor--Wiles method, one usually assumes the existence of a suitable~$\m$ as a hypothesis.
The usual shorthand way of describing what comes out of the Taylor--Wiles method is then an~``$R = \T$ theorem.'' 
Proving an~$R =\T$ theorem  can more or less be divided into three different problems:

\begin{enumerate} 
\item Understanding~$\T$. Why does there exist\footnote{At the time of Wiles' result, this was seen as the easier
direction (if not easy), although, in light of the success of the Taylor--Wiles method, it
may well be the harder direction in general.}
 a map~$R \rightarrow \T$?
This is the problem of the ``existence of Galois representations.''  Implicit here is the problem of showing that those Galois representations not only exist but have
the ``right local properties''  at the ramified primes, particularly those dividing~$p$.
\item Understanding~$R$. Wiles introduced a  mechanism for controlling~$R$ via its
tangent space using Galois cohomology (in particular Poitou--Tate duality~\cite{Milne}),
and this  idea  has  proved remarkably versatile. 
What has changed, however, is our understanding of local Galois representations and how this
information can be leveraged to understand the structure of~$R$.
\item Understanding why the  map~$R \rightarrow \T$ is an isomorphism. 
\end{enumerate}

We begin by summarizing the original~$R=\T$ theorem from this 
viewpoint (or more precisely, the modification by Faltings which appears as an appendix to~\cite{TW}). We only discuss for now the so-called
 ``minimal case\footnote{The case when  the Galois representations attached to~$R$ and~$\T$  have minimal level~$N$ as determined by  the residual representation.}'' since this is most relevant for subsequent generalizations (see~\S\ref{section:kisin}).
Our summary is cursory, but see~\cite{MR1605752,CSS}  for  excellent expositional sources on  early versions of the Taylor--Wiles method.
We start with a representation~$\rhobar: G_{\Q} \rightarrow \GL_2(\Fbar_p)$ for~$p > 2$ which (say) comes from a semistable elliptic curve~$E$
and which we assume to be modular.
By a theorem of Ribet~\cite{Ribet100}, we may assume it is modular of level either~$N = N(\rhobar)$ or~$N = N(\rhobar)p$ where~$N(\rhobar)$ is the Serre
weight~\cite{MR885783} of~$\rhobar$.

\begin{enumerate}
\item Understanding~$\T$: The construction of Galois representations associated to modular forms has its own interesting history (omitted here), 
but (in the form originally needed by Wiles)  was
more or less complete for modular forms (and even Hilbert modular forms)  by 1990. The required local properties at primes different from~$p$ followed from work of Carayol~\cite{Carayol},
and the local properties at~$p$ were well understood either by Fontaine--Laffaille theory~\cite{FL}, or, in the ordinary case,
by Mazur--Wiles~\cite{MW} (see also work of Hida~\cite{MR848685,MR1463699}).
\item Understanding~$R$:  Here~$R$ is a deformation ring of~$\rhobar$ subject to precise local deformation conditions at~$p$
and the primes dividing~$N(\rhobar)$. For the prime~$p$,
the local conditions amount either to  an ``ordinary'' or ``finite-flat'' restriction.
One then interprets the dual of the reduced tangent space~$\m_R/(\m^2_R,p)$ of~$R$ in terms of Galois cohomology,
in particular  as a subgroup (Selmer group) of classes in~$H^1(\Q,\ad^0(\rhobar))$ satisfying  local conditions. This can be thought
of as analogous to a class group, and one does not have any \emph{a priori} understanding of how large it can be  although it has some finite dimension~$d$.
Using the Greenberg--Wiles formula, the obstructions  in~$H^2(\Q,\ad^0(\rhobar))$ can
be related to the reduced tangent space, and  allow one to realize~$R$ as a quotient
of~$W(k) \llbracket x_1, \ldots, x_d \rrbracket$ by $d$ relations. In particular, if~$R$ 
was finite and free as a~$W(k)$-module (as would be the case if~$R = \T$) then~$R$ would be a complete intersection.
\item Understanding why the  map~$R \rightarrow \T$ is an isomorphism. Here lies the heart of the Taylor--Wiles method.
The ring~$\T$ acts on a natural module~$M$ of modular forms. One shows --- under a mild hypothesis on~$\rhobar$  
--- the existence of (infinitely many) sets~$Q = Q_N$  for any natural number~$N$ of cardinality~$|Q| = d$  --- so-called Taylor--Wiles primes --- with a number of pleasant properties:
\begin{enumerate} 
\item The primes~$q \in Q$ are congruent to~$1 \bmod p^N$.
\item Let~$R_Q$ be the deformation ring capturing the same local properties as~$R$ but modified so that the representations at primes in~$Q$
may now be ramified of degree~$p^N$.
There is naturally a surjection~$R_{Q} \rightarrow R$, but for Taylor--Wiles primes,  this modification does not increase the size of the tangent space.
In particular, for a fixed ring~$R_{\infty} = W(k) \llbracket x_1, \ldots, x_d \rrbracket$ there are surjections~$R_{\infty} \rightarrow R_Q \rightarrow R$ for every~$Q$.
\item The corresponding rings~$\T$ and~$\T_{Q}$ act naturally on spaces of modular forms~$M$ and~$M_Q$ respectively. 
Using multiplicity one theorems, Wiles proves (see~\cite[Theorem~2.1]{W})  that~$M$ and~$M_Q$ are free of rank one over~$\T$ and~$\T_Q$ respectively. The space~$M$ can be interpreted as a space of  modular forms for
a particular modular curve~$X$. The second key property of Taylor--Wiles primes is that there are no new modular forms associated to~$\rhobar$ at level~$X_0(Q)$, and hence~$M$ can \emph{also}
be interpreted as a  space of  modular forms for~$X_0(Q)$. 
There is a Galois cover~$X_1(Q) \rightarrow X_0(Q)$ with Galois group~$(\Z/Q \Z)^{\times}$, and hence an intermediate cover~$X_H(Q) \rightarrow X_0(Q)$
with Galois group~$\Delta_N = (\Z/p^N \Z)^d$ acting via diamond operators. The space~$M_Q$ is  essentially a localization of a certain space of modular forms for~$X_H(Q)$
(with some care taken at the Hecke operators for primes dividing~$Q$).
Since the cohomology of modular curves (localized at the maximal ideal corresponding to~$\m$)
is concentrated in degree one, the module~$M_Q$ turns out to be free over an auxiliary ring~$S_N = W(k)[\Delta_N]$ of diamond operators,
and the quotient~$M_Q/\a_Q$ for the augmentation ideal~$\a_Q$ of~$S_N$ is isomorphic to~$M$. It follows that~$\T_Q/\a_Q = \T$.
\item The diamond operators have an interpretation on the Galois deformation side, and there is a identification~$R_Q/\a_Q = R$ where~$R_Q$ and~$\T_Q$ can be viewed compatibly as~$S_{N}$-modules.
\end{enumerate}
\item Finally, one ``patches'' these constructions together for larger and larger~$Q$. This is somewhat counterintuitive, since for different~$Q$
the Galois representations involved are not compatible. However, one forgets the Galois representations and only remembers the structures relative
to both the diamond operators~$S_N$ and~$R_{\infty}$,  giving the data of a surjection
$$R_{\infty} \rightarrow \T_{\infty}$$
with a  compatible action of~$S_{\infty} = \projlim S_N \simeq W(k)\llbracket t_1, \ldots, t_d \rrbracket$.
Using the fact that~$\T_{\infty}$ is free of finite rank over~$S_{\infty}$,
and that~$R_{\infty}$ and~$S_{\infty}$ are formally smooth of the same dimension, one deduces that~$R_{\infty}
=\T_{\infty}$ and then~$R = \T$ after quotienting out by the augmentation ideal of~$S_{\infty}$.
\end{enumerate}

\section{The early years}

\subsection{The work of Diamond and Fujiwara}

Wiles 
made essential uses of multiplicity one theorems in order to deduce that~$M_Q$ was free over~$\T_{Q}$.
Diamond~\cite{DiamondMult} and Fujiwara~\cite{Fujiwara} (independently)
had the key insight that one could instead patch the modules~$M_Q$ directly --- and then argue directly with the resulting object~$M_{\infty}$ as a module over~$R_{\infty}$ which was also
free over~$S_{\infty}$.
Using the fact that~$R_{\infty}$ is formally smooth, this allowed
one to deduce \emph{a posteriori} that~$M_{\infty}$ was free over~$R_{\infty}$
using the Auslander--Buchsbaum 
formula~\cite{ABF}. This not only removed the necessity of
proving difficult multiplicity one results but gave new proofs of these 
results\footnote{There is an intriguing result of Brochard~\cite{Brochard} which weakens the hypotheses of Diamond's freeness criterion even further,
although this idea has not yet been fully exploited.}
 which could then be generalized to situations where the known methods
(often using the~$q$-expansion principle) were unavailable\footnote{The history of the subject involves difficult theorems in the arithmetic geometry
of Shimura varieties being replaced by insights from commutative algebra, paving the way to generalizations where further insights from the arithmetic
geometry of Shimura varieties are required.}. Diamond had the following to say about how he came up with the idea to patch
 modules rather than use multiplicity one theorems:
\smallskip

\begin{footnotesize}
My vague memory is that I was writing down examples of ring homomorphisms and modules, subject to some constraints imposed by a Taylor--Wiles setup, and I couldn't break  ``$M$ free over the group ring implies~$M$ free over~$R$.''  (I still have the notebook with the calculations somewhere, mostly done during a short trip with some friends to Portugal.)  I didn't know  \emph{what} commutative algebra statement I needed, but I knew I needed to learn more commutative algebra and found my way to Bruns and Herzog's ``Cohen-Macaulay Rings'' \cite{MR1251956} (back in the library in Cambridge UK by then).  When I saw the statement of Auslander--Buchsbaum, it just clicked.

\end{footnotesize}

\smallskip
Diamond made a second  improvement~\cite{MR1405946,MR1638490} dealing with primes away from~$p$ in situations where the corresponding minimal
local deformation problem was not controlled by the Serre level~$N(\rhobar)$ alone.

\subsection{Integral \texorpdfstring{$p$}{p}-adic Hodge Theory, part I: Conrad--Diamond--Taylor} \label{CDT}
One early goal after Fermat was the resolution of the full Taniyama--Shimura conjecture, namely,
the modularity of all elliptic curves over~$\Q$.
After the improvements of Diamond, the key remaining problem was  understanding 
deformation rings associated to local Galois representations at~$p$  coming from elliptic curves with bad reduction at~$p$.
Since Wiles' method (via Langlands--Tunnell~\cite{MR574808,MR621884})
was ultimately reliant on working with the prime~$p = 3$, this meant understanding deformations at~$p$ of level~$p^2$
and level~$p^3$, since any elliptic curve over~$\Q$ has a twist such that the largest power of~$3$ dividing the conductor is at most~$27$.
Ramakrishna in his  thesis~\cite{RamakrishnaThesis} had studied the local deformation problem for finite flat representations (the case when~$(N,p) = 1$) and
proved that the corresponding local deformation rings were formally smooth.
The case when~$p$ exactly divides~$N$ was subsumed into the ordinary case, also treated by Wiles.
In level~$p^2$,  one can show that the Galois representations associated to the relevant
 modular forms\footnote{This is not true for all modular forms of level~$p^2$ and weight~$2$,
but only for those whose conductor at~$p$ remains divisible by~$p^2$ after any quadratic twist.}
  of level~$p^2$ become finite flat
after 
passing to a finite extension~$L/\Q_p$ with ramification degree~$e \le p- 1$.
In this range, Conrad~\cite{Conrad1,Conrad2} was able to adapt ideas of Fontaine~\cite{Fontaine} to give an equivalence between the local Galois
deformations (assuming~$\rhobar |_{G_{\Q_p}}$ was irreducible) and linear algebra data. In particular, 
as in the work of Ramakrishna, one can show that the relevant local deformation rings are formally smooth, and so
Conrad, Diamond, and Taylor were able to adapt
the Taylor--Wiles method  to this setting~\cite{CDT}.

\subsection{Integral \texorpdfstring{$p$}{p}-adic Hodge Theory, part II: Breuil--Conrad--Diamond--Taylor} \label{BCDT}

A central technical ingredient in all of the arguments so far has been some use of integral $p$-adic Hodge Theory,
and in particular the theory of finite flat group schemes and Barsotti--Tate groups as developed by Fontaine and others.
All  integral versions of this theory required a hypothesis on either the weight or the ramification index~$e$ relative to the bound~$p-1$.
However, around this time, Christophe Breuil made a 
breakthrough\footnote{Much of the development of integral~$p$-adic Hodge theory over the last~$20$ years since~\cite{BCDT}
has been inspired by  its use in the Taylor--Wiles method. However, the timing of Breuil's work was more of a happy coincidence,
although Breuil was certainly aware of the fact that a computable theory of finite flat group schemes over highly ramified bases could well have
implications in the Langlands program.}
 by finding a new way to understand the integral theory
of finite flat group schemes over arbitrarily ramified bases~\cite{BreuilAnnals}. This was just the technical tool required
to push the methods of~\cite{CDT} to level~$p^3$. 
Using these results, Breuil, Conrad, Diamond and Taylor~\cite{BCDT} were able to show that  enough suitably chosen local deformation
rings were formally smooth to prove the modularity of all elliptic curves.

\subsection{Higher weights, totally real fields, and base change}

Many of the methods which worked for modular forms were directly adaptable to the case both of general rank~$2$
motives over~$\Q$ with distinct Hodge--Tate weights (corresponding to modular forms of weight~$k \ge 2$ rather than~$k = 2$)
and also to such motives over totally real fields (which are related to Hilbert modular forms), see
in particular the work of Fujiwara~\cite{Fujiwara} (and more recently Frietas--Le Hung--Siksek~\cite{BAOBAO}).
Another very useful innovation was a base change idea of Skinner--Wiles~\cite{SWBC}
which circumvented the need to rely on Ribet's level 
 lowering theorem. The use of cyclic base change (\cite{MR574808} in this case and~\cite{MR1007299} in general) subsequently became a standard tool in the subject.
 For example, it  meant that one could always reduce to a situation where the ramification at all primes~$v \nmid p$ was unipotent.
The paper~\cite{SWBC} was 
 related to a more ambitious plan  by Wiles to prove modularity for all totally real fields:

\smallskip

\begin{footnotesize}
After Fermat I started to work with Taylor and then Diamond on the general case but decided very soon that I would rather try to do the totally real case for~$\GL(2)$.
I think this was while I was getting back into other kinds of problems but I thought I should still earn my bread and butter. One lunch time at the IAS in 1996 Florian Pop spoke to me and explained to me about finding points over fields totally split at some primes (e.g. real places) as he had written a paper~\cite{POP} about this with some others.
 Was this any use for the Tate--Shavarevich group? I immediately saw that whether or not it was any use for TS (I doubted it) it should certainly give potential modularity. This gave some kinds of lifting so I worked on the other half (i.e. descent) thinking that just needed a similar insight. At some point I suggested to Chris that we try to do Ribet's theorem using cyclic base change as that would be useful in proving modularity and was buying time while I waited to get the right idea. Unfortunately I completely misjudged the difficulty of descent and the problem is still there. I think it is both much harder than I thought and also more important. I hope still to prove it! Of course Taylor found potential modularity and then, what I had assumed was much harder, a way to think about~$\GL(n)$.
 
 \end{footnotesize}

\section{Reducible representations: Skinner--Wiles}
\label{section:SW}
One of the key  hypotheses in the Taylor--Wiles method concerns restrictions
on the representation~$\rhobar$, in particular the hypothesis that~$\rhobar |_{G_{\Q(\zeta_p)}}$
is absolutely irreducible. In~\cite{SW1,SW2,SW3}, Skinner and Wiles introduced
a new argument in which this hypothesis was relaxed, at least assuming
the representations were ordinary at~$p$. In the ordinary setting, one can replace the
rings~$R$ and~$\T$ (which in the original setting are finite over~$W(k)$) by rings which
are finite  (and typically flat) over Iwasawa algebras~$\Lambda = W(k) \llbracket (\Z_p)^d  \rrbracket$ for some~$d$
which arise as weight spaces,
the point being that the ordinary deformations of varying weight admit a good integral theory.
The first innovation (in part) involves making a base
change so that the reducible locus is (relatively) ``small,'' (measured in terms of the
codimension over~$\Lambda$). The second idea is then to apply a variant of the Taylor--Wiles
method to representations~$\varrho: G_F \rightarrow \GL_2(\T/\mathfrak{p})$
for non-maximal prime ideals~$\mathfrak{p}$ of~$F$\footnote{Representations~$\varrho$
 to infinite quotients~$\T/\pp$ had also arisen in Wiles'  paper on
Galois representations associated to ordinary modular forms~\cite{WilesOrdinary} 
 where the concept of pseudo-deformation was also first introduced.}. Wiles again:
 
 \smallskip
 
 \begin{footnotesize}
 We had worked out a few cases we could do without big Hecke rings in some other papers and I would say it was more a feat of stamina and technique to work through it. Of course the use of these primes was much more general and systematic than anything that went before. There is also an amusing point in this paper where we use a result from commutative algebra. It seemed crucial then though I don't know if it still is. This is proposition A.1 of Raynaud~\cite{MR0407021}. I had thought at some point during the work on Fermat that this result might be needed and had asked Michel Raynaud about it. He said he would think about it. A week later he came back to me, somewhat embarrassed that he had not known right away, to say that it was a result in his wife's thesis. So the reference to M.Raynaud is actually to his wife, Mich\`{e}le Raynaud,  though he gave the reference.
 
 \end{footnotesize}

 Allen~\cite{Allen} was later able to adapt these arguments to  the~$p=2$ dihedral case, which
 (in a certain sense)  realized the original 
desire\footnote{As far as primary historical sources go, the introduction of Wiles' paper~\cite{W}
is certainly worth reading.}  of Wiles to work at the prime~$p=2$.

\section{The Artin Conjecture}
\label{section:artin}
While the approach of~\cite{W,TW} applied  (in principle) to all Galois representations
associated to modular forms of weight~$k \ge 2$, the case of modular forms of weight~\mbox{$k = 1$}
is qualitatively quite different (see also~\S\ref{section:calegari}). 
It was therefore quite surprising when Buzzard--Taylor~\cite{MR1709306}
proved weight one modularity lifting theorems for  odd continuous representations~$\rho: G_{\Q}
\rightarrow \GL_2(\Qbar_p)$ which were unramified at~$p$.
Using this, Buzzard--Dickinson--Shepherd-Barron--Taylor~\cite{MR1845181}  proved the Artin conjecture for a positive proportion of all odd~$A_5$
representations, which had previously only been known in a finite number of cases\footnote{In a computational \emph{tour de force} for the time, Buhler~\cite{Buhler} in his thesis
had previously established the modularity of an explicit odd projective~$A_5$ representation of conductor~$800$.} up to twist.
Standard ordinary modularity theorems showed the existence of ordinary modular forms
associated to such representations~$\rho$ --- however, the classicality theorems of
Hida~\cite{MR848685} do not apply (and are not true!) in weight one.
The main idea of~\cite{MR1709306} was to exploit the fact that~$\rho$ is unramified to construct
\emph{two} ordinary modular forms each corresponding to a choice
of eigenvalue of~$\rho(\Frob_p)$ assuming these eigenvalues are 
distinct\footnote{This argument can be modified to deal with the case when
the eigenvalues of~$\rhobar(\Frob_p)$ coincide by  modifying~$R$ and~$\T$
to include operators corresponding (on the Hecke side) to~$U_p$.
Geraghty and I discovered an integral version of this  idea ourselves (``doubling,'' following Wiese's paper~\cite{Wiese})
during the process of writing~\cite{CG}, although it turned out that, at least in characteristic zero, Taylor already
had the idea in his back pocket in the early 2000s.}.
One then  has to argue~\cite{MR1709306} that these two ordinary forms are the oldforms associated
to a classical eigenform of weight one, which one can do by exploiting
both the rigid geometry of modular curves and the~$q$-expansion principle. 

Although the original version of this argument required a number of improvements
to the usual Taylor--Wiles method (Dickinson overcame some technical issues  when~$p=2$~\cite{MR1845182} and
Shepherd-Barron--Taylor proved some new cases of Serre's conjecture for~$\SL_2(\F_4)$ and~$\SL_2(\F_5)$-representations
 in~\cite{MR1491977}), it was ripe
for generalization to totally real fields\footnote{The proof all that finite odd~$2$-dimensional representations over~$\Q$ are modular
was completed by Khare and Wintenberger
as a consequence of their proof of Serre's conjecture, see~\S\ref{sec:Serre}.}
After a key early improvement by
Kassaei~\cite{MR2219265}, the~$n = 2$ Artin conjecture for totally real fields is now
completely resolved under the additional assumption that the representation is 
odd by a number of authors, including Kassaei--Sasaki--Tian and Pilloni--Stroh~\cite{MR3581174,MR2983010,MR3294624,SasakiInvent,MR3639600,MR3581178}.
On the other hand, the reliance on~$q$-expansions in this 
argument
 has proved
an obstruction to extending this to other 
groups.
(See also~\S\ref{section:GSp4}).

\section{Potential Modularity}
\label{section:potential} \label{section:pq}
One new idea which emerged in Taylor's paper~\cite{MR1954941} was the concept of \emph{potential modularity}. Starting with a representation~$\rho: G_F \rightarrow \GL_2(\Qbar_p)$ for a totally
real field~$F$,
one could sidestep the (difficult) problem of proving the modularity of~$\rhobar$ by proving
it was modular over some finite totally real extension~$F'/F$.   In the original paper~\cite{W},
Wiles employed a~$3$-$5$ switch to deduce the modularity of certain mod~$5$ representations
from the modularity of mod~$3$ representations. More generally,
one can prove the modularity of a mod~$p$ representation~$\rhobar_p$
from the modularity of a mod~$q$ representation~$\rhobar_q$ if one can find both of them
occurring as the residual representation of a compatible family where the Taylor--Wiles hypotheses
apply to~$\rhobar_q$. For example, if~$\rhobar_p$ and~$\rhobar_q$ are 
representations valued in~$\GL_2(\F_p)$ and~$\GL_2(\F_q)$ respectively, one can try to
find the compatible family by finding an elliptic curve with a given mod~$p$ and mod~$q$ representation.
The obstruction to doing such  a~$p$-$q$ switch over~$F$
is
that the corresponding moduli spaces (which in this case are twists of the modular curve~$X(pq)$)  are not in general rational, and hence have no
reason to admit rational points. However, exploiting an idea
due to Moret-Bailly~\cite{MB}, Taylor showed that these moduli spaces at least had many points over totally
real fields where one could additionally ensure that the Taylor--Wiles hypothesis applies at the prime~$q$.
At the cost of proving a weaker result, this gives a huge amount of extra flexibility that has proved
remarkably useful. Taylor's first application of this idea was to prove the Fontaine--Mazur
conjecture for many~$2$-dimensional representations, since the potential modularity of these representations
 was enough to prove (for example) that they come from compatible families of Galois representations
 (even over the original field~$F$!), and that they satisfy purity (which is  known for Hilbert modular
 forms of regular weight). The concept of potential modularity, however, has proved crucial for other applications, not least of which is the proof of the Sato--Tate conjecture (see~\S\ref{section:satotate}).

\section{The work of Kisin}

A key ingredient in the work of Breuil--Conrad--Diamond--Taylor (\S\ref{CDT},\S\ref{BCDT})
(and subsequent work of Savitt~\cite{MR2004122,Savitt}) 
was the fact that a certain local deformation ring~$\Rfl$ defined in terms of integral~$p$-adic
Hodge theory was formally smooth. The calculations of~\cite{BCDT,Savitt}, however, applied only to some
(very) carefully chosen situations sufficient for elliptic curves but certainly not for all~$2$-dimensional
representations. In the 2000s, Kisin made a number of significant contributions, both to the understanding
of local deformation rings but also to the structure of the Taylor--Wiles argument itself~\cite{Kisin1,Kisin2,Kisin3,Kisin4,Kisin5,Kisin6}.

\subsection{Local deformation rings at \texorpdfstring{$v = p$}{v eq p}}
One difficulty with understanding local deformation rings~$\Rfl$ associated to finite flat group schemes
over  highly ramified  bases is that the group schemes themselves are not uniquely defined by their generic fibres.
Kisin~\cite{Kisin4} had the idea that  one could also define the moduli space of the group schemes
themselves, giving a projective resolution~$\RGR \rightarrow \mathrm{Spec}(\Rfl)$ (this map is
an isomorphism after inverting~$p$).
Kisin further realized that the geometry of~$\RGR$ was related to local models
of Shimura varieties, for which one had other available techniques to analyze their structure and singularities.
Later, Kisin was also able~\cite{Kisin6} to construct  local deformation rings~$R$  capturing deformations
of a fixed local representation~$\rhobar$ which become  semi-stable over a fixed extension~$L/\Q_p$
and had Hodge--Tate weights in any fixed finite range~$[a,b]$, absent  a complete integral theory of  such representations.
(There are also are constructions where one fixes the inertial type of the corresponding representation.)
Kisin further proved that the generic fibres of these rings were indeed of the expected dimension and often 
 formally smooth.
 
\subsection{Kisin's modification of Taylor--Wiles} \label{section:kisin}
Beyond analyzing the local deformation rings themselves, Kisin crucially found a way~\cite{Kisin4} to modify the Taylor--Wiles method to
avoid the requirement that these rings are formally smooth, thus greatly expanding the scope of the method.
First of all, Kisin reimagined the global deformation ring~$R$ as an algebra
over a (completed tensor product) 
\begin{equation*}
\Rloc = \widehat{\bigotimes}_{v \in S} R_v
\end{equation*}
of local deformation rings~$R_v$ for sets of places~$v \in S$, in particular  including the
 prime~$p$\footnote{Since the local residual representations are typically reducible,
 Kisin also introduced the notion of \emph{framed} deformation rings which are always well-defined,
 and which (properly taking into account the extra variables) are compatible with the Taylor--Wiles
 argument.}.
Now,  after a Taylor--Wiles
patching argument, one constructs a big module~$M_{\infty}$ over~$R_{\infty}$
(and free over the auxiliary ring of diamond operators~$S_{\infty}$)
but where~$R_{\infty}$ is no longer a power series ring over~$W(k)$
but a power series ring over~$\Rloc$.  If the algebras~$R_v$ for~$v \in S$ are themselves power series rings, one is reduced
precisely to the original Taylor--Wiles setting as modified by Diamond. On the other hand, if the~$R_v$ are (for example)
not power series rings but are
integral domains over~$W(k)$ of the expected dimension,
then Kisin explained how one could still deduce that~$M[1/p]$ was a faithful~$R[1/p]$-module, which proves
that~$R[1/p]=\T[1/p]$ and suffices for applications to modularity. 
More generally, assuming only that the~$R_v$ are  flat over~$W(k)$ and that the generic fibre~$R_v[1/p]$ is 
equidimensional of  the expected dimension,  the modularity of any point of~$R$ 
reduces to showing that there is at least one modular point which
lies on the same component of~$R_v[1/p]$\footnote{There are some subtleties to understanding~$R[1/p]$ for complete local Noetherian~$W(k)$-algebras
that are not obvious on first consideration. The first and most obvious blunder to avoid is to recognize that~$R[1/p]$
is usually far from being a local ring. Similarly, the ring~$R[1/p]$ can be regular
 and still have multiple components,
as can be seen in an example as simple as~$R=\Z_p\llbracket X \rrbracket/X(X-p)$.
Perhaps more importantly, however, the ring~$R[1/p]$ ``behaves''  in some important ways like a finitely generated algebra over a field.}.

In the original modularity lifting arguments, one treated the minimal case first and then deduced the non-minimal cases
using a subtle commutative algebra criterion which detected isomorphisms between complete intersections.
From the perspective of Kisin's modification, all that is required is to show that there exists a single modular point
with the right non-minimal local properties. In either case, both Wiles and Kisin used Ihara's Lemma
 to establish the
existence of congruences between old and new forms, but Kisin's argument is much softer and thus more generalizable
to other situations\footnote{In particular, Wiles' numerical criterion~\cite[Thm~5.3]{MR1605752} relies on certain rings being complete
intersections, and Kisin's local deformation rings are not complete intersections (or even Gorenstein) in general --- see~\cite{Snowden}.}.
Kisin had the following to say about his thought process:

\smallskip

\begin{footnotesize}
The idea of thinking of~$R$ as an~$\Rloc$ algebra just popped into my head, after I'd been thinking about the Wiles--Poitou--Tate formula, and how it fit into the Taylor--Wiles patching argument. This was in Germany, I think in 2002. I had the idea about moduli of finite flat group schemes in the Fall of 2003, after I arrived in Chicago. 
It was entirely motivated by modularity. I had been trying to compute these deformation rings, by looking at deformations of finite flat group schemes. For~$e<p-1$, the finite flat model is unique, so I knew this gave the deformation ring in this case; this already gave some new cases. However I was stuck about the meaning of these calculations in general for quite some time. At some point I thought I'd better write up what I had, but as soon as I started thinking about that --- within a day --- I realized what the correct picture was with the families of finite flat group schemes resolving the deformation ring. I already knew about Breuil's unpublished note~\cite{Breuilnote}, and quite quickly was able to prove the picture was correct. It was remarkable that prior to coming to Chicago, I didn't even know the definition of the affine Grassmannian, but within a few months of arriving, it actually showed up in my own work.

To me the whole project was incredibly instructive. If I had known more about what was (thought to be) essential in the Taylor--Wiles method, I never would have started the project. Not having fixed ideas gave me time to build up intuition. I also should have gotten the idea about moduli of finite flat group schemes much sooner if I'd been more attentive to what the geometry was trying to tell me.

\end{footnotesize}

\vspace{-1em}

\section{\texorpdfstring{$p$}{p}-adic local Langlands}

\subsection{The Breuil--M\'{e}zard conjecture}

Prior to Kisin's work,  Breuil and M\'{e}zard~\cite{BM} undertook a study of certain low weight potentially semi-stable
deformation rings, motivated by~\cite{BCDT}. They discovered (in part conjecturally)  a crucial link between the geometry of these Galois deformation rings
(in particular, the Hilbert--Samuel multiplicities of their special fibres) with  the mod-$p$ reductions (and corresponding
irreducible constituents) of lattices inside locally algebraic~$p$-adic representations of~$\GL_2(\Z_p)$. 
In the subsequent papers~\cite{Breuil1,Breuil2}, Breuil 
 raised the hope that there could exist a~$p$-adic Langlands correspondence relating certain mod-$p$ (or~$p$-adic Banach space)
representations of~$\GL_2(\Q_p)$
to geometric~$2$-dimensional~$p$-adic representations of~$G_{\Q_p}$\footnote{The starting observation~\cite{BreuilIntro} is as follows: if~$\pi = \bigotimes' \pi_v$
is the automorphic representation associated to a modular form~$f$, then~$\pi_v$ determines (and is determined by)~$\rho_f  |_{G_{\Q_v}}$ for all~$v \ne p$ (at least
up to Frobenius semi-simplification).  On  the other hand,  $\pi_p$ does not determine the~$p$-adic representation~$\rho_f |_{G_{\Q_p}}$
(except in the exceptional setting where~$\pi_p$ is spherical and~$a_p$ is not a~$p$-adic unit),  raising the question of what extra~$\GL_2(\Q_p)$ structure associated
to~$f$ should
determine (and be determined by)~$\rho_f |_{G_{\Q_p}}$.}. Breuil recounts the origins of these conjectures as follows:

\medskip

{\footnotesize

The precise moment I became 100\% sure that there would be a non-trivial~$p$-adic correspondence for~$\GL_2(\Q_p)$ was in the computations of~\cite{Breuil2}. In these computations, I reduced mod~$p$ certain~$\Zbar_p$--lattices in certain locally algebraic representations of~$\GL_2(\Q_p)$, and at some point, I found out that this reduction mod~$p$ had a really nice behaviour, so nice that clearly, it was predicting (via the mod-$p$ correspondence) what the reduction mod-$p$ would be on the~$\Gal(\Qbar_p/\Q_p)$-side.

}

These ideas were further developed by Colmez in~\cite{MR2493219,MR2642407,MR2642409} amongst other papers\footnote{In~\cite{KisinEarly}, Kisin had shown
 that the~$p$-adic representations~$V$ associated to non-classical finite slope overconvergent modular forms with~$U_p$-eigenvalue~$a_p$ satisfied~$\dim D_{\mathrm{cris}}(V) = 1$, and moreover that crystalline
Frobenius acted on this space by~$a_p$.
(This
paper was itself apparently motived by the goal of \emph{disproving} the Fontaine--Mazur conjecture!)
On the way to the~$2004$ Durham symposia on~$L$-functions and Galois representations, Fontaine
raised the question to Colmez to
 what extent this determined the corresponding Galois representation.
Colmez worked out the answer the evening before his talk 
and incorporated it into his lecture the following day,  ultimately leading to the notion of trianguline representations~\cite{MR2493219}.}:
Colmez studied  various Banach space completions defined by Breuil and proved they
were non-zero using the theory of~$(\varphi,\Gamma)$-modules. 
 Since the theory of~$(\varphi,\Gamma)$-modules applies to all Galois representations and not 
just potentially semi-stable ones, this led Colmez to propose a~$p$-adic local Langlands correspondence for
\emph{arbitrary}~$2$-dimensional representations~$G_{\Q_p} \rightarrow \GL_2(E)$, and he was ultimately able to construct a functor
from suitable~$\GL_2(\Q_p)$-representations to Galois representations of~$G_{\Q_p}$.
Colmez gave a talk on his construction at a conference in Montreal in September 2005. At the same conference, Kisin gave a  talk
presenting a proof of the Breuil--M\'{e}zard conjecture by relating it directly to~$R=\T$ theorems and the Fontaine--Mazur conjecture for
 odd~$2$-dimensional
representations of~$G_{\Q}$ with distinct Hodge--Tate weights. While Kisin's argument exploited results of Berger--Breuil~\cite{BB} and Colmez, it was realized  
by the key participants (perhaps in real time) that Colmez' $p$-adic local Langlands correspondence should be viewed as taking place over the entire local deformation ring.
Subsequently Colmez was able to construct the inverse functor\footnote{To add some further confusion to the historical chain of events,
the published version of~\cite{KisinFM} incorporates some of these subsequent developments. Note also that the
current state of affairs is that the proof of the full~$p$-adic local Langlands correspondence for~$\GL_2(\Q_p)$ (for example as proved in~\cite{MR3272011}
but see also~\cite[Remarque VI.6.51]{MR2642409})
still relies on the global methods of~\cite{emerton2010local}, which in turn relies on~\cite{MR2642409}.
These mutual dependencies, however, are not circular! The difficulty arises in the supercuspidal case. One philosophical reason that global methods
are useful here is that all global representations are yoked together by an object (the completed cohomology group~$\widetilde{H}^1(\Z_p)$) with
good finiteness properties. One can then exploit  the fact
 that crystabeline  representations (for which the~$p$-adic local Langlands correspondence is known by~\cite{MR2642409}) 
 are Zariski dense inside unrestricted global deformation rings
(\cite[Theorem~1.2.3]{emerton2010local},  using arguments going back to B\"{o}ckle~\cite{MR1854117}).}.
Colmez writes:

\medskip

{\footnotesize

 I received a paper of Breuil (a former version of~\cite{MR2667887}) during my stay at the Tata Institute in December 2003---January 2004. In December, 
 I was spending Christmas under Goa's palm trees with my daughter when Breuil's paper arrived in my email. That paper contained a conjecture (in the semi-stable case) that I was sure I could prove using~$(\varphi,\Gamma)$-modules (if it was true\ldots). I spent January 2004 working on it and after 15 days of computations in the dark, I finally found a meaning to some part of a painful formula (you can find some shadow of all of this in (iii) of Remark 0.5 of my unpublished~\cite{ColmezUnpublished}). By the end of the month, I was confident that the conjecture was proved and I told so to Breuil who adapted the computations to the crystalline case, and wrote them down with the help of Berger (which developed into~\cite{BB}). (One thing that makes computations easier and more conceptual in the crystalline case is that you end up with the universal completion of the locally algebraic representation you start with; something that is crucial in Matthew [Emerton]'s proof of the FM conjecture.)  Durham was in August of that year and Berger--Breuil had notes from a course they had given in China~\cite{BBnotes}.
 Those notes were instrumental in my dealing with trianguline representations at Durham (actually, I did some small computation and the theory just developed by itself during the night before my talk which was supposed to be on something else\ldots I think I came up with the concept of trianguline representations later, to justify the computations, I don't remember what language I used in my talk which had some part on Banach--Colmez spaces as far as I can remember.
 
 }

 \subsection{Local--global compatibility for completed cohomology}
 
 From a different perspective, Emerton had introduced the completed cohomology groups~\cite{MR2207783}
 as an alternative means for constructing the Coleman--Mazur eigencurve~\cite{CM}. Inspired by Breuil's work, Emerton
  formulated~\cite{EmertonOne}
 a local--global compatibility conjecture for completed cohomology in the language of the then nascent~$p$-adic Langlands correspondence.
 After the construction of the correspondence for~$\GL_2(\Q_p)$ by Colmez and Kisin, Emerton was able to prove most of his conjecture, leading to a new proof of (many cases of) the  Fontaine--Mazur conjecture.
 The results of Kisin~\cite{KisinFM} and Emerton fell short of proving the full version of this conjecture for two reasons. The first was related to some technical 
 issues with the~$p$-adic local Langlands correspondence, both at the primes~$p=2$ and~$3$ but also when the residual representation 
 locally had the shape~$1 \oplus \overline{\varepsilon}$ for the cyclotomic character~$\varepsilon$. (The local issues have now more
 or less all been resolved~\cite{MR3272011}. The most general  global results for~$p=2$ are currently due to Tung~\cite{MR4257080}.)
A second restriction was the Taylor--Wiles
 hypothesis that~$\rhobar$ was irreducible. Over the intervening years, a number of key improvements to the local story have been found,
 in particular by Colmez, Dospinescu, Hu, and Pa{\v{s}}k{\=u}nas~\cite{MR3150248,MR3272011,PH}.
 Very recently, Lue Pan~\cite{Pan1} found a way to marry techniques from Skinner--Wiles in the reducible case (\S\ref{section:SW})  to techniques from~$p$-adic local Langlands to
 completely prove the modularity (up to twist) of any geometric representation~$\rho: G_{\Q} \rightarrow \GL_2(\Qbar_p)$ for~$p \ge 5$ only assuming the hypotheses
 that~$\rho$ has distinct Hodge--Tate weights and that~$\rho$ is odd\footnote{The assumption on the Hodge--Tate weights is almost certainly removable 
 using recent progress on the ideas  discussed in~\S\ref{section:artin} (Sasaki has announced such a result).  
 Moreover, Pan has found a different
 approach to this case as well,  see~\cite[Theorem 1.0.5]{Pan2} and the subsequent comments.
 The hypothesis that~$\rho$ is odd is more troublesome --- see~\S\ref{section:even}.}.

\section{Serre's Conjecture} \label{sec:Serre}
In Wiles' original lectures in Cambridge in 1993, he introduced 
his method with the statement that it was~\emph{orthogonal} to Serre's conjecture~\cite{MR1322785}.
 In some senses, this viewpoint turned out to be the opposite of prophetic,
in that the ultimate resolution of Serre's conjecture used the Taylor--Wiles  method as its  central core.
The proof of Serre's conjecture by Khare and Wintenberger~\cite{SerreKW1,SerreKW2,KW1,KW2} 
introduced a new technique
for lifting residual Galois representations to characteristic zero (see~\S\ref{section:KW}) which
has proved very useful for subsequent modularity lifting theorems.

\subsection{Ramakrishna lifting}
Ramakrishna, in a series of papers in the late 90s~\cite{RK1,RK2}, studied the question of lifting an odd Galois representation
\begin{equation*}
\rhobar:G_{\Q} \rightarrow \GL_2(\Fbar_p)
\end{equation*}
to a global potentially semistable representation in characteristic zero  unramified outside finitely many primes. This is a trivial consequence
of Serre's conjecture\footnote{Trivial only assuming the results of Tsuji~\cite{Tsuji} and Saito~\cite{Saito}, of course.}
but is highly non-obvious without such an assumption. Ramakrishna succeeded in proving the existence of lifts by an ingenious
argument involving adding auxiliary primes and modifying the local deformation problem to a setting where there all global obstructions
vanished. The resulting lifts had the added property that they were valued in~$\GL_2(W(k))$ whenever~$\rhobar$ was valued in~$\GL_2(k)$.
Adaptations of Ramakrishna's method had a number of important applications even under the assumption of residual modularity, including in~\cite{CHT} where it was used to produce
characteristic zero lifts with Steinberg conditions at some auxiliary primes. 
There is also recent work of Fakhruddin, Khare, and Patrikis~\cite{FKP} which considerably
extends these results in a number of directions.

\subsection{The Khare--Wintenberger method} \label{section:KW}
One disadvantage of Ramakrishna's method was that it required allowing auxiliary ramification which (assuming Serre's
conjecture) should not be necessary\footnote{If one insists on finding a lift is valued in~$\GL_2(W(k))$ rather than~$\GL_2(\OL_E)$
for some ramified extension~$E/W(k)[1/p]$, then some auxiliary ramification is necessary in general, at least in fixed weight.}. 
Khare and Wintenberger found a new powerful method for avoiding this. The starting point is the idea
that, given an odd representation~$\rhobar: G_F \rightarrow \GL_2(\Fbar_p)$ for a totally real field~$F$ satisfying the Taylor--Wiles hypotheses,
one could find a finite extension~$H/F$ where~$\rhobar$ is modular (exactly as in~\S\ref{section:potential}). Then, using an~$R=\T$ theorem over~$H$,
one proves that the corresponding deformation ring~$R_H$ of~$\rhobar |_{G_H}$ is finite over~$W(k)$. However, for formal reasons, there is a map~$R_H \rightarrow R_F$
(where~$R$ is the deformation ring corresponding to the original representation~$\rhobar$) which is a finite morphism, and hence the ring~$R_F/p$ is Artinian.
Then, by Galois cohomological arguments, one proves the ring~$R_F$ has dimension at least one, from which one deduces that~$R_F$ has~$\Qbar_p$-valued points.
Even more can be extracted from this argument, however --- the~$\Qbar_p$-valued point of~$R_F$ certainly comes from a~$\Qbar_p$-valued point of~$R_H$,
and hence comes from a compatible family of Galois representations over~$H$. Using the fact that one member of the family extends to~$G_F$, it can be argued that
the entire family descends to a compatible family over~$F$. This one can then hope to prove is modular by working at a different (possibly smaller) prime, where (hopefully)
one can prove the associated residual representation is modular. In this way, one can inductively
reduce Serre's conjecture~\cite{MR885783} to the case~$p=2$ and~$N(\rhobar)=1$,  where Tate had previously proved in a letter to Serre~\cite[July 2, 1973]{ColmezTate2} (also~\cite{Tate}) that all such absolutely irreducible 
representations are modular by showing that no such representations exist. The entire idea is very clean, although in practice the difficulty reduces to the step
of proving modularity lifting theorems knowing either that~$\rhobar$ is either modular and absolutely irreducible or is reducible. Khare and Wintenberger's
timing was such that the automorphy lifting technology was just good enough for the proof to work, although this required some extra effort at the prime~$p = 2$
(both in their own work and in a key assist by Kisin~\cite{Kisin3}).
As with Ramakrishna's method, the Khare--Wintenberger lifting method has also been systematically exploited for modularity lifting applications
(for example in~\cite{BLGGT} (see~\S\ref{section:BLGGT}) building on ideas of Gee~\cite{gee061}).

\section{Higher dimensions}

Parallel to the developments of~$p$-adic Langlands for~$n = 2$, the first steps were made to generalize the theory to higher dimensional
representations. 
Unlike in the case of modular forms,  substantially less was known about the  existence of Galois representations until the 90s.

\subsection{Construction of Galois representations, part I: Clozel--Kottwitz}

The first general construction\footnote{Clozel's paper is from 1991 and thus not strictly ``post--Fermat'' as is the remit of this survey.
However, it can be considered a natural starting point for the ``modern'' arithmetic theory of automorphic forms for~$\GL(n)$ and so
it seems
 reasonable to mention it here.}
of Galois representations in dimension~$n > 2$ was made by Clozel~\cite{Clozel} (see also the work of Kottwitz~\cite{Kottwitz}).
Clozel's theorem
 applies to certain automorphic forms for~$\GL_n(\A_L)$ for CM fields~$L/L^{+}$.
The construction requires
three important hypotheses on~$\pi$: The first is 
that~$\pi$ is conjugate self-dual, that is~$\pi^{\vee} \simeq \pi^c$. If~$\pi$ is a base change from an algebraic representation of~$L^{+}$ and~$n = 2$ then
this condition is automatic\footnote{At least after a twist which is always possible to achieve in practice,
see~\cite[Lemma~4.1.4]{CHT}.
More generally, one can work with unitary similitude groups and consider~$\pi$ with~$\pi^{\vee} \simeq \pi^c \otimes \chi$
for suitable characters~$\chi$.},
 but it is far from automatic when~$n > 2$. The second condition is an assumption
on the infinitesimal character which (in the case of modular forms) is equivalent to the condition that the weight~$k$ is~$\ge 2$.
Finally, there is a technical condition that for some finite place~$x$ the representation~$\pi_x$
is square integrable. A number of improvements (particularly at the bad primes) were made by Harris--Taylor
in~\cite[Theorem~C]{ht}, and later by Taylor--Yoshida and Caraiani~\cite{TY,MR2972460,MR3272276}, bringing the theory roughly in line with that of modular forms at the time of Wiles,
and in particular primed for possible generalizations of the Taylor--Wiles method to higher dimensions.

\subsection{The Sato--Tate conjecture, part I}
\label{section:satotate}

Harris, and Taylor  (as early as 1996) started the work of generalizing the Taylor--Wiles machinery to the setting of~$n$-dimensional representations.
They quickly understood that the natural generalization of these ideas in~$n$-dimensions required the hypothesis
that the Galois representations were self-dual up to a twist. 
This meant that one should not consider general automorphic forms on the group~$\GL_n(\A_{\Q})$ but rather groups of symplectic or orthogonal type
depending on the parity of~$n$. If one replaced Galois representations over totally real fields by Galois representations over imaginary CM
fields and then further imposed the condition that the Galois representations are conjugate self-dual, the relevant automorphic forms should then come
from unitary groups.
There were two benefits of working with these hypotheses. First of all, the relevant automorphic representations for unitary groups were, as with modular forms,
associated to cohomology classes on Shimura varieties. In particular, under the assumption that there existed
an auxiliary prime~$x$ such that~$\pi_x$ was square integrable, they could
be seen inside the ``simple'' Shimura varieties of type~$U(n-1,1)$ considered by Kottwitz~\cite{Kottwitz}. On the other hand, the same
Hecke eigenclasses (if not Galois representations) also came from a compact form of the group and thus inside the cohomology
of zero-dimensional varieties\footnote{Inside~$H^0$, of course}. The advantage of working in this setting is that the freeness of~$M_Q$
over the ring of diamond operators is immediate\footnote{In more general contexts, the freeness of~$M_Q$ is closely related to the vanishing of cohomology
localized at~$\m$ in all but one degree.}.
In the fundamental paper~\cite{CHT}, Clozel, Harris, and Taylor succeeded in overcoming many of the technical difficulties generalizing the arguments
of~\cite{W,TW} to these representations. Although the argument in spirit was very much the same, there are a number of points for~$\GL_2$
where things are much easier.  
One representative example of this  phenomenon is understanding Taylor--Wiles primes. While the Galois side generalizes readily, the automorphic side requires many new ideas
and some quite subtle arguments concerning the mod~$p$ structure of certain~$\GL_n(\Q_q)$-representations of conductor~$1$ and conductor~$q$.
In order to prove the Sato--Tate conjecture for a modular form~$f$, it was already observed by Langlands that it sufficed to prove the modularity
of all the symmetric powers of~$f$. However, it turns out that the weaker assumption that each of these symmetric powers is \emph{potentially modular}
suffices, and by some subterfuge only the even powers are required~\cite{Harristwist}. In order to prove potential modularity theorems, one
needs to be able to carry out some version of the~$p$-$q$ switch (\S\ref{section:pq}). In order to do \emph{this}, one needs a source
of motives which both generate Galois representations of the right shape (conjugate self dual and with distinct Hodge--Tate weights)  and yet also come in positive dimensional families.
It turned out that  there already existed such motives in the literature, namely,  the so-called Dwork family. However, given the strength of the automorphy
lifting theorems in~\cite{CHT}, considerable effort had to be made in studying the geometry of the Dwork family to ensure that the~$p$-$q$ switch
would produce geometric Galois representations with the right local properties. These issues were precisely addressed in the companion
paper by Harris, Shepherd-Barron and Taylor~\cite{Shep}. Taken together, these papers contained all the ingredients to prove the potential modularity of higher symmetric
powers of modular forms (satisfying a technical square integrable condition at some auxiliary prime) with one exception.
As mentioned earlier, the work of Kisin had simplified the passage from the minimal case to the non-minimal case --- ``all'' that was required
was to produce congruences between the original form and forms of higher level rather than to compute a precise congruence number as in~\cite{W}.
However, even applying Kisin's approach seemed to require Ihara's Lemma, and despite several years of 
effort, the authors of~\cite{CHT} were not able to 
overcome this obstacle\footnote{The issue remains unresolved to this day.}. Here is Michael Harris' recollection of the process:

\medskip

{\footnotesize

	In the spring of 1995, I was at Brandeis, Richard was at MIT, and I wanted 
to understand the brand new proof of Fermat's Last Theorem.   So I asked Richard
if he would help me learn by collaborating on modularity for higher-dimensional groups.   The collaboration
took off a year later, when Richard wrote to tell me about the Diamond--Fujiwara
argument and suggested that we work out the Taylor--Wiles method for unitary groups.  This
developed over the next 18 months or so into the early version of what eventually
became the IHES paper with Clozel.  But it had no punch line. 
 I was hoping to work out
some non-trivial examples of tensor product functoriality for~$\GL(n)\times \GL(m)$, where one of the 
two representations was congruent mod~$l$ to one induced from a CM Hecke character.  
This would have required some numerical verification.  In the meantime we got sidetracked into 
proving the local Langlands conjecture~\cite{ht}.  

	The manuscript on automorphy lifting went through several drafts and was circulated;
you can still read it on my home page~\cite{Harrissecret}.  Genestier and Tilouine~\cite{GT} quoted it when they proved 
modularity lifting for Siegel modular forms.  When Clozel saw the draft he told me we should try to prove 
the Sato--Tate Conjecture.   Although this was in line with my hope for examples of tensor product functoriality, 
it seemed completely out of reach, because I saw no way to prove residual modularity of symmetric powers.

When I heard about the Skinner--Wiles paper I came up with a quixotic plan to prove symmetric power 
functoriality for Eisenstein representations, using the main conjecture of Iwasawa theory to
control the growth of the deformation rings.  This was in the spring of 2000, at the IHP special 
semester on the Langlands program, where I first met Chris Skinner. 

	One day Chris told me that Richard had invented potential modularity.  This led me to a 
slightly less hopeless plan to prove potential symmetric power functoriality by proving it for~$2$-dimensional 
representations congruent to potentially abelian representations, as in the potential modularity 
argument.  I told Richard about this idea, probably the day he arrived in Paris.  He asked:  why
apply potential modularity to the~$2$-dimensional representation; why not instead apply it to
the the symmetric power representations directly?  I then replied:  that would require a variation
of Hodge structures with a short list of properties:  mainly, the correct~$h^{p,q}$'s and large monodromy groups.
We checked that potential modularity was sufficient for Sato--Tate.  We then resolved to ask our 
contacts if they knew of VHS with the required properties.  The whole conversation lasted about 20 minutes.

	I asked a well-known algebraic geometer, who said he did not know of any such VHS.  
Richard asked Shepherd-Barron, who immediately told him about the Calabi-Yau hypersurfaces 
that had played such an important role in the mirror symmetry program.  
(And if my algebraic
geometer hadn't wanted to be dismissive, for whatever reason, he would have realized this as well.)  
The~$h^{p,q}$'s were fine but we didn't know about the monodromy.  However, Richard was staying at the 
IHES, and by a happy accident so was Katz, and when Richard asked Katz about the monodromy 
for this family of hypersurfaces Katz told him they were called the Dwork family and gave him the 
page numbers in one of his books.

So within a week or two of our first conversation, we found ourselves needing only one more
result to complete the proof of Sato--Tate.  This was Ihara's lemma, which occupied our attention over 
the next five years.  In the meantime, Clozel had written a manuscript on symmetric powers,
based on the reducibility mod ell of symmetric powers.  The argument was incomplete but he had 
several ideas that led to his joining the project, and he also hoped to use ergodic theory to prove Ihara's lemma.
In the summer of 2003 Clozel and I joined Richard in ``old'' Cambridge to try to work this out.  The rest
you know.   We finally released a proof conditional on Ihara's lemma in the fall of 2005.  A few months
later Richard found his local deformation argument, and the proof was complete.

}

\subsection{Taylor's trick: Ihara Avoidance}
\label{section:iharaavoid}
Shortly after the preprints~\cite{CHT,MR1491977} appeared, Taylor found a way to overcome the problem
of Ihara's lemma. Inspired  by Kisin's formulation of the Taylor--Wiles method (\S\ref{section:kisin}), Taylor
had the idea of comparing two global deformation rings~$R^1$ and~$R^2$. Here (for simplicity) the local deformation
problems associated to~$R^1$ and~$R^2$ are formally smooth at all but a single prime~$q$. At the prime~$q$, however,
the local deformation problem associated to~$R^1$ consists of tamely ramified representations where a generator~$\sigma$
of tame inertia has  characteristic polynomial~$(X-1)^n$, and for~$R^2$ the
characteristic polynomial has the shape~$(X-\zeta_1)\ldots (X - \zeta_n)$
for some fixed distinct roots of unity~$\zeta_i \equiv 1 \bmod p$. On the automorphic side, there are two patched modules~$H_1$ and~$H_2$,
and there is an equality~$H_1/p = H_2/p$. The local deformation ring~$R^1_q$ associated to~$R^1$ at~$q$ is reducible
and has multiple components
in the generic fibre, although the components in characteristic zero are in bijection to the components in the special fibre.
 On the other hand, the local deformation ring associated to~$R^2$ at~$q$ consists of a single component, and so using Kisin's
 argument one deduces that~$H_2$ has full support. Now a commutative algebra argument using the identity~$H_1/p=H_2/p$
 and the structure of~$R^1_q$ implies that~$H_1$ has sufficiently large support over~$R_1$, enabling one  to deduce the modularity
 of every~$\Qbar_p$-valued point of~$R_1$\footnote{Taylor's argument proves theorems of the form~$R[1/p]^{\mathrm{red}} = \T[1/p]$ rather than~$R = \T$. This
is still perfectly sufficient for proving modularity lifting results, but not always other interesting corollaries associated to~$R=\T$ theorems
 like finiteness of the corresponding adjoint Selmer groups (though see~\cite{PDuke,1912.11265}).}.

\subsection{The Sato--Tate conjecture, part II}

After Taylor's trick, one was \emph{almost} in a position to complete the proof of Sato--Tate for all classical modular forms.
A few more arguments were required. One was the tensor product trick due to Harris which enabled one to pass
from conjugate self-dual motives with weights in an arithmetic progression to conjugate self-dual motives with consecutive
Hodge--Tate weights by a judicious twisting argument using CM characters. A second ingredient was the analysis of
the ordinary deformation ring by Geraghty~\cite{ger}.  One of the requirements of the~$p$-$q$ trick was the  condition that
 certain moduli spaces (the Dwork family in this case) had points over various local extensions~$E$ of~$\Q_p$, in order to construct
 a motive~$M$ over a number field~$F$ with~$F_v = E$ for~$v|p$.  For the purposes of modularity
lifting, one wants strong control over the local deformation ring at~$p$, and the choice of local deformation ring is more or less forced by
the geometric properties of the~$p$-adic representations associated to~$M$ . One way to achieve this would be to work in
the Fontaine--Laffaille range where the local deformation rings were smooth. But this requires both that~$M$ is smooth at~$p$
and that the ramification degree~$e$
of~$E/\Q_p$ is one. It is not so clear, however, that the Dwork family contains suitable points (for a fixed residual representation~$\rhobar$)
which lie in any unramified extension of~$\Q_p$. What Geraghty showed, however, was that certain ordinary deformation rings\footnote{In Geraghty's setting,
the residual representations~$\rhobar$ were locally trivial. Hence the definition of ``ordinary''  was not something that could be defined
on the level of Artinian rings, and the construction (as with Kisin's construction of local deformation rings associated to certain types) is therefore indirect.}
  were connected over arbitrarily ramified bases. 
The final piece, however, was the construction of Galois
representations for all conjugate self-dual regular algebraic cuspidal~$\pi$ without the extra condition that~$\pi_q$
was square integrable for some~$q$. This story merits it own separate discussion; suffice to say
that it required the combined efforts of many people and the resolution of many difficult problems, not least of which
was  the fundamental lemma by Laumon and Ng\^{o}~\cite{LN,N} (see also the Paris book project~\cite{MR2742611,book2},
the work of Shin~\cite{Shin},
and many more references which if I attempted to make complete
would weigh down the bibliography and still contain grievous omissions).

\subsection{Big image conditions}

The original arguments in~\cite{W,TW} required a ``big image'' hypothesis, namely that~$\rhobar$ was absolutely irreducible after restriction to the Galois group of~$\Q(\zeta_p)$.
Wiles' argument also required the vanishing of certain cohomology groups associated to the adjoint representation of the image of~$\rhobar$.
These assumptions had  natural analogues in~\cite{CHT} (so-called ``big image'' hypotheses) although they were quite restrictive, and it wasn't clear that they would
even apply to most residual representations coming from some irreducible compatible family.
In the setting of~$2$-dimensional representations, the Taylor--Wiles
hypothesis guarantees the existence of many primes~$q$ such that~$q \equiv 1 \bmod p$ and such that~$\rhobar(\Frob_q)$ has distinct eigenvalues.
This ensures, for example, that there cannot be any Steinberg deformations at~$q$ because the ratio of the eigenvalues of any Steinberg deformation must be~$q$.
In dimension~$n$, one natural way to generalize this might be to say that~$\rhobar(\Frob_q)$ has distinct eigenvalues, although this is not always possible
to achieve for many irreducible representations~$\rhobar$. A weaker condition is that~$\rhobar(\Frob_q)$ has
an eigenvalue~$\alpha$ with multiplicity one. For such~$q$,  there will be no deformations which are unipotent on inertia at~$q$ for which the generalized~$\alpha$
eigenspace is not associated to a~$1$-dimensional block. The translation of this into an automorphic condition on~$U_q$-eigenvalues  is precisely what is done in~\cite{CHT}
(there are additional technical conditions on~$\Frob_q$ with respect to the adjoint representation~$\ad(\rhobar)$ which we omit here).
In~\cite{jack}, however, Thorne finds a way to allow~$\rhobar(\Frob_q)$ to have an eigenvalue~$\alpha$ with  higher multiplicity, and yet still cut out (integrally) the space
of automorphic forms whose Galois representations decompose at~$q$ as an unramified representation plus a one dimensional representation which is tamely ramified of~$p$-power
order. This technical improvement is very important because (as proved in the appendix by Guralnick, Herzig, Taylor, and Thorne~\cite{jack}) it imposes no restrictions on~$\rhobar$ when~$p \ge 2n+1$ 
beyond the condition that~$\rhobar$ is
absolutely irreducible after restriction to~$G_{\Q(\zeta_p)}$. This improvement is very useful for applications.

\subsection{Potentially diagonalizable representations} \label{section:BLGGT}

After the proof of Sato--Tate for modular forms, Barnet-Lamb, Gee, and  Geraghty turned their attention to proving the analogous theorem for Hilbert modular forms of regular weight.
The methods developed so far were well-suited both to representations~$\rho$ which were either ordinary or when~$\rho$ was not ordinary but still Fontaine--Laffaille.
(The latter implies that~$\rhobar |_{G_{\Q_p}}$ is absolutely irreducible  of some particular shape.)
For a modular form over~$\Q$, one easily sees that~$\rho$  takes one of one of these forms for any sufficiently large~$p$. 
For Hilbert modular forms, one certainly expects that the ordinary hypothesis should hold for all~$v|p$ and infinitely many~$p$, but this remains open. The difficulty arises when, for some prime~$p$
(that splits completely, say)
the~$p$-adic representation is ordinary at some~$v|p$ but
non-ordinary other~$v|p$.
The reason that this causes issues
is that, when applying the Moret-Bailly argument in the~$p$-$q$ switch, one wants to avoid any ramification at~$p$ for the non-ordinary case, and yet have large ramification
at the ordinary case to make~$\rhobar$ locally trivial, and these desires are not compatible. The resolution in~\cite{HilbertST} involved a clever refinement of the Harris
tensor product trick. These ideas were further refined in~\cite{BLGGT} and led to the concept of a \emph{potentially diagonalizable} 
representation~$\rho: G_E \rightarrow \GL_n(\Qbar_p)$ for some finite extension of~$E/\Q_p$. Recall from~\S\ref{section:kisin} that, in the modified form of the Taylor--Wiles
method, proving modularity of some lift of~$\rhobar$ often comes down to showing the existence of a modular lift lying on a smooth point of the corresponding component of
the generic fibre of~$\Rloc$. In light of Taylor's Ihara avoidance trick (\S\ref{section:iharaavoid}), the difficulty in this  problem is mostly at the prime~$p$, and in particular
the fact that one knows very little about the components of general Kisin potentially crystalline deformation rings. A potentially diagonalizable representation is one for which,
after some finite (necessarily solvable!) extension~$E'/E$, the representation~$\rho |_{G_E'}$ is crystalline and lies on the same generic irreducible component as a diagonal representation.
This notion has a number of felicitous properties. First, it includes Fontaine--Laffaille representations and ordinary potentially crystalline representations. Second, it is clearly invariant under base change.
Third, it is compatible with the tensor product trick of Harris. These features make it supremely well-adapted to the current forms of the Taylor--Wiles method.
By combining this notion with methods of~\cite{HilbertST,blght}, as well as extensive use of Khare--Wintenberger lifting (\S\ref{section:KW}), Barnet--Lamb, Gee, Geraghty, and Taylor in~\cite{BLGGT}
proved the potential automorphy of \emph{all} conjugate self-dual irreducible\footnote{One variant proved shortly thereafter by Patrikis--Taylor~\cite{PT} replaced
the irreducibility condition by a purity condition (which is automatically satisfied by representations coming from pure motives).} odd\footnote{Although
there is no longer a non-trivial complex conjugation in the Galois group of a CM field, there is still an oddness condition related to the conjugate self-duality of the representation
and the fact that there are two ways for an irreducible representation to be self-dual (orthogonal and symplectic).}
 compatible systems of Galois representations over a totally real field.

\subsection{Even Galois Representations} 
\label{section:even}
The Fontaine--Mazur conjecture for geometric Galois representations~$\rho: G_{\Q} \rightarrow \GL_2(\Qbar_p)$ predicts that, up to twist, either~$\rho$
is modular \emph{or} $\rho$ is even with finite image. The methods of~\cite{W,TW} required the assumption that~$\rhobar$ was modular
and so \emph{a priori} the assumption that~$\rho$ was odd (at least when~$p > 2$).
Nothing at all was known about the even case before the papers~\cite{Cal1,Cal2} in which a very simple trick made the problem accessible
to modularity lifting machinery under the assumption that the Hodge--Tate weights are distinct. The punch line is that, for any CM field~$F/F^{+}$,
the restriction~$\Sym^2(\rho): G_{F} \rightarrow \GL_3(\Qbar_p)$ is conjugate self-dual and 
no longer sees the ``evenness'' of~$\rho$\footnote{The representation~$\rho$ itself restricted to~$G_F$ will not be odd in the required sense --- 
one exploits the fact here that~$3$ is odd whereas symplectic representations are always even dimensional.}. Hence one can hope to prove it is potentially modular for some CM extension~$L/L^{+}$,
and then by cyclic base change~\cite{MR1007299} potentially modular for the totally real field~$L^{+}$. But  Galois representations coming from
regular algebraic automorphic forms for totally real fields will not be 
even\footnote{I managed to twist Taylor's arm into writing the paper~\cite{TaylorSign} which proved this for odd~$n$, which
sufficed for my purposes where~$n$ was either~$3$ or~$9$. This is now also known for general~$n$, see Caraiani--Le Hung~\cite{AnaBao}.},
 and thus one obtains a contradiction. 
 These ideas are already enough to deduce the main result of~\cite{Cal1} directly from~\cite{BLGGT}, although in contrast~\cite{Cal2} uses (indirectly) the full strength of the~$p$-adic local Langlands correspondence
via theorems of Kisin~\cite{KisinFM}. The papers~\cite{Cal1,Cal2} still fall short of completely resolving the Fontaine--Mazur in this case even for~$p > 7$,
since there remain big image hypotheses on~$\rhobar$.
On the other hand, this trick has nothing to say about the case when the Hodge--Tate weights are equal (see~\S\ref{section:hopeless}).

\subsection{Modularity of higher symmetric powers}

Another parallel development in higher dimensions was the extension of Skinner--Wiles (\S\ref{section:SW}) to higher dimensions.
Many of the arguments of Skinner--Wiles relied heavily on the fact that any proper submodule of a~$2$-dimensional representation must
have dimension~$1$, and one-dimensional representations are very well understood by class field theory.
Nonetheless, in~\cite{MR3327536}, Thorne proved a residually reducible modularity theorem for higher dimensional
representations. In order to overcome the difficulty of controlling
reducible deformations, he imposed a  Steinberg condition at some auxiliary prime. Although this is a definite restriction,
it does apply (for example) to the  Galois representation coming from the symmetric power of a modular form which also satisfies this condition.
In a sequence of papers~\cite{CT1,CT2,CT3}, Clozel and Thorne applied this modularity lifting theorem to prove new cases of 
symmetric power functoriality (see also the paper of Dieulefait~\cite{Dieulefait}). A key difficulty here is again the absence of Ihara's lemma in order to find automorphic forms
with the correct local properties. Very recently (using a number of new ideas), Newton and Thorne~\cite{NT1,NT2} were able to (spectacularly!) complete this program and prove the
full modularity of all symmetric powers of all modular forms.

\section{Beyond self-duality and Shimura varieties}

All the results discussed so far --- with the exception of those discussed in~\S\ref{section:artin} ---- apply only to Galois representations which 
are both regular and satisfy some form of self-duality. Moreover, they all correspond to automorphic forms which can be detected by the (\'{e}tale) cohomology
of Shimura varieties. Once one goes beyond these representations, many of the established methods begin to break
 immediately\footnote{I should warn the reader that this section and the next (even more than the rest of this paper) is filtered through the lens of my own
 personal research journey --- \emph{caveat lector}!}.

An instructive case to consider is the case of~$2$-dimensional geometric Galois representations of an imaginary quadratic field~$F$ with distinct Hodge--Tate weights.
The corresponding automorphic forms for~$\GL_2(\A_F)$ contribute to the cohomology of locally symmetric spaces~$X$ which are arithmetic 
hyperbolic~$3$-manifolds\footnote{Already by 1970, Serre (following ideas of Langlands) was trying to link Mennicke's computation that~$\GL_2(\Z[\sqrt{-109}])^{\mathrm{ab}}$ is infinite to the possible existence of elliptic curves over~$\Q(\sqrt{-109})$ with good reduction everywhere~\cite[Jan 14, 1970]{ColmezTate1}.}.
These spaces are certainly not algebraic varieties and their cohomology is hard to access via algebraic methods.
One of the first new questions to arise in this context is the relationship between torsion classes and Galois representations.
Some speculations about this matter were made by Elstrodt, Grunewald, and Mennicke at least as far back as 1981~\cite{EGM}, but the most influential conjecture was due to Ash~\cite{Ash},
who conjectured that eigenclasses in the cohomology of congruence subgroups of~$\GL_n(\Z)$ over~$\Fbar_p$ (which need not lift to characteristic zero)
 should give rise to~$n$-dimensional Galois 
representations over finite fields.
Later, conjectures were made~\cite{AshSerre,AshExtra} in the converse conjecture in the spirit of Serre~\cite{MR885783} linking Galois representations to classes in cohomology modulo~$p$.
Certainly around~2004, however, it was not at all clear what exactly one should expect the landscape to be\footnote{I recall conversations  with a number of experts at the 2004 Durham conference, where nobody seemed quite sure even what the dimension of the ordinary deformation ring~$R$ of a~$3$-dimensional representation~$\rho: G_{\Q} \rightarrow \GL_3(\Fbar_p)$ should be. Ash, Pollack, and Stevens had computed numerical examples where a regular algebraic ordinary cuspidal form for~$\GL_3(\A_{\Q})$ 
not twist-equivalent to a symmetric square did not appear to admit classical deformations. (I learnt about this example from Stevens at a talk at Banff in
December~$2003$.)
This would be easily explained if~$R$ had (relative) dimension~$0$ over~$\Z_p$
but be more mysterious otherwise.}, and so it was around this time I decided to start thinking about this question\footnote{One great benefit  to me at the time of thinking about
Galois representations over imaginary quadratic fields was that it didn't require me to understand the geometry of Shimura varieties which I have always found
too complicated to understand. The irony of course is that the results of~\cite{10author,BCGP} ultimately rely on extremely intricate properties of Shimura varieties.}
 in earnest.  I became convinced very soon (for aesthetic reasons if not anything else) that if one modified~$\T$ to be the
ring of endomorphisms acting on integral cohomology (so that it would see not only the relevant automorphic forms but also the torsion classes) then there should still
be an isomorphism~$R=\T$. Moreover, this equality would not only be a form of reciprocity which moved beyond the conjecture linking motives to automorphic forms,
but it suggested that the integral cohomology of arithmetic groups (including the torsion classes) were themselves the fundamental object of interest.
Various developments served only to confirm this point of view. In my paper with Mazur~\cite{CalMazur}, we gave some theoretical evidence for why ordinary
families of Galois representations of imaginary quadratic fields might on the one hand be positive dimensional and explained completely by torsion classes and
yet not contain any classical automorphic points at all. During the process of writing~\cite{CalDunfield}, Dunfield (numerically) compared the torsion classes
in the cohomology of inner forms of~$\GL_2$ and the data was in perfect agreement with a conjectural Jacquet--Langlands correspondence for torsion
(later taken up in joint work with Venkatesh~\cite{CAVE}). Emerton and I had the idea of working with completed cohomology groups both to construct Galois representations and
even possibly to approach questions of modularity.  The first idea was to exploit the well-known relationship between the cohomology of these manifolds
and the cohomology of the boundary of certain Shimura varieties. We realized that if we could control the co-dimension of the completed cohomology
groups over the non-commutative Iwasawa algebra, the Hecke eigenclasses would be forced to be seen by eigenclasses coming from the middle
degree of these Shimura varieties where one had access to Galois representations\footnote{Unfortunately, these conjectures~\cite[Conj~1.5]{MR2905536} remain all open in more
or less all cases except for Scholze's results in the case of certain Shimura varieties~\cite[Cor~4.2.3]{scholze-torsion}.}. 
On the automorphy lifting side,
we had even vaguer ideas~\cite[\S1.8]{MR2905536}\footnote{Pan's remarkable paper~\cite{Pan1}  turned some of these pipe dreams
into reality.} on how to proceed. A different (and similarly unsuccessful) approach\footnote{This is taken from my 2006 NSF proposal, and I believe
influenced by my conversations with Taylor at Harvard around that time.}  was to work with ordinary deformations over a partial weight space for a split prime~$p=vw$ in an imaginary quadratic field~$F$.
That is, deformations of~$\rhobar$ which had an unramified quotient at~$v$ and~$w$ but with varying weight at~$v$ and fixed weight at~$w$.
Here the yoga of Galois deformations suggested that~$R$ should be finite flat over~$W(k)$ in this case (and even a complete intersection). Moreover,
one had access to~$\T$ using an overlooked\footnote{One should never overlook results of Hida. I only learnt about this paper when Hida pointed it out to me
(with a characteristic smile on his face)
after my talk in Montreal in 2005. I was pleased at least that the idea that these families were genuinely non-classical was not anticipated either
in~\cite{HidaIQF} or  in~\S4 of Taylor's thesis~\cite{taylorthesis}.}
 result of Hida~\cite{HidaIQF},  and in particular one could deduce that~$\T$ has dimension at least one.
If one could show that~$\T$ was flat over~$W(k)$, then one could plausibly apply (assuming the existence of Galois representations) the original argument of~\cite{W,TW}.
The flatness of~$\T$, however, remains an open problem\footnote{One might even argue that there is no compelling argument to believe it is true ---
the problem is analogous to the vanishing of the~$\mu$-invariant in Iwasawa theoretic settings.}.

\subsection{The Taylor--Wiles method when \texorpdfstring{$l_0 > 0$}{l g 0}, part I: Calegari--Geraghty}
\label{section:calegari}

Shortly before (and then during) the special year on Galois representations at the IAS in 2010-2011, I started to work with Geraghty in earnest on the problem
of proving~$R=\T$ in the case of imaginary quadratic fields, assuming the existence of a surjection~$R \rightarrow \T$.
A computation in Galois cohomology shows that the expected ``virtual'' dimension of~$R$ over~$W(k)$ should be~$-1$, and hence the patched module~$M_{\infty}$
should have codimension~$1$ over the ring of diamond operators~$S_{\infty}$. 
 We realized this was a consequence of the fact that, after localizing the cohomology
at a non-Eisenstein maximal ideal, the cohomology should be non-zero in exactly two degrees. More precisely, patching the presentations
of these~$S_N$-modules
 would result in a balanced presentation of~$M_{\infty}$ as an~$S_{\infty}$-module with the same  (finite) number
of generators as relations.
We then
realized that the same principle held more generally for~$n$-dimensional representations over any number field. In characteristic zero, the 
localized cohomology
groups were non-zero exactly in a range~$[q_0,q_0+l_0]$ (with~$q_0$ and~$l_0$ as defined in~\cite{BM}) where~$-l_0$ coincided with the expected virtual dimension of~$R$ over~$W(k)$ coming
from Galois cohomology. We could thus show --- assuming the localized torsion cohomology also vanished in this range --- that  by
patching \emph{complexes} $P_Q$ (rather than modules~$M_Q$),
one arrives at a complex~$P_{\infty}$ of free~$S_{\infty}$ modules in degrees~$[q_0,q_0 + l_0]$. 
Because the ring~$R_{\infty}$ of dimension~$\dim R_{\infty} = (\dim S_{\infty}) - l_0$ acts  by patching on~$H^*(P_{\infty})$, a
 simple commutative
algebra lemma then shows that~$M_{\infty} = H^*(P_{\infty})$ has codimension~$l_0$ over~$S_{\infty}$
and must be concentrated in the final degree. In particular, the Taylor--Wiles method (as modified by Diamond) could be happily adapted 
to this general setting\footnote{David Hansen came up with a number of these ideas independently~\cite{Hansen}.}.
Moreover, the arguments were compatible with all the other improvements, including Taylor's Ihara avoidance argument amongst other things\footnote{These methods
only prove~$R[1/p]^{\mathrm{red}} = \T[1/p]^{\mathrm{red}}$, of course. In situations where~$\T \otimes \Q = 0$, the methods of~\cite{CG} in the minimal
case prove not only that~$R=\T$ but also that (both) rings are complete intersections. Moreover, one also has access to level raising (on the level of complexes)
and Ihara's lemma~\cite[\S4]{CAVE}, and I tried for some time (unsuccessfully) to adapt the original minimal $\Rightarrow$ non-minimal arguments
of~\cite{W} to this setting. There certainly seems to be some rich ideas in commutative algebra in these situations to explore, see for example
recent work of Tilouine--Urban~\cite{TU}.}.
We also realized that the same idea applied to Galois representations coming from the coherent cohomology  of Shimura varieties even when
the corresponding automorphic forms were not discrete series. While our general formulation involved a number of conjectures we considered hopeless,
the coherent case had at least one setting in which many more results were available, namely the case of modular forms of weight one,
where the required vanishing conjecture was obvious, and where we were able to establish the existence of the required map~$R \rightarrow \T$
with all the required local properties by direct arguments. Although the state of knowledge concerning Galois representations
increased tremendously between the original conception of~\cite{CG} and its final publication, by early 2016 it still seemed out of reach
to make any of the results in~\cite{CG} unconditional.

\subsection{Construction of Galois representations, part II}

Before one can hope to prove~$R=\T$ theorems, one needs to be able to associate Galois representations to the corresponding automorphic forms.
There are two contexts in which one might hope to make progress. The first is in situations where the automorphic forms contribute to the 
Betti cohomology of some locally symmetric space --- for example, tempered algebraic cuspidal automorphic representations for~$\GL_n(\A_F)$
and any~$F$. The second is in situations  where the automorphic forms contribute to the coherent cohomology of some Shimura variety.
Here the first and easiest case corresponds to weight one modular forms, where the Galois representations were first constructed
by Deligne and Serre~\cite{deligne-serre}.

In work of Harris--Soudry--Taylor~\cite{HTS,HTST}, Galois representations were constructed for regular algebraic forms for~$\GL_2(\A_F)$ for an imaginary quadratic field~$F$ and satisfying a further
restriction on the central character. Harris, Soudry, and Taylor exploited (more or less) the fact that the automorphic induction of such forms 
are self-dual (although not regular) and still contribute to coherent cohomology,
so one can construct Galois representations using a congruence argument as in the paper of Deligne and Serre~\cite{deligne-serre}.
On the other hand, this does not prove the expected local properties of the Galois representation at~$v|p$.

It was well-known for many years that the Hecke eigenclasses associated to regular algebraic cuspidal automorphic forms for~$\GL_n(\A_F)$
for a CM field~$F$  could be realized as eigenclasses coming from the boundary of certain unitary Shimura varieties of type~$U(n,n)$. It was, however,
also well-known that the corresponding \'{e}tale cohomology classes did not realize the desired Galois 
representations\footnote{For a more basic example of what can go wrong, note that the Hecke eigenvalues of~$T_{l}$
on~$H^0(X,\Q_p)$ of a modular curve are~$1+l$, which corresponds to the Galois representation~$\Q_p \oplus \Q_p(1)$. However, only the piece~$\Q_p$
occurs inside~$H^0$.}.  Remarkably, this problem was completely and unexpectedly resolved in 2011 in~\cite{HLTT} by Harris--Lan--Taylor--Thorne.
Richard Taylor writes:

\smallskip

\begin{footnotesize}
For~\cite{HLTT} I knew that the Hecke eigenvalues we were interested in contributed to Betti cohomology of~$U(n,n)$. The problem was to show that they contributed to overconvergent $p$-adic cusp forms. I was convinced on the basis of Coleman's paper ``classical and overconvergent modular forms'' \cite{MR1369416} that this must be so. I can't now reconstruct exactly why Coleman's paper convinced me of this, and it is possible, even probable, that my reasoning didn't really make any sense. However it was definitely this that kept me working at the problem, when we weren't really getting anywhere.

\end{footnotesize}

Amazingly, this breakthrough  immediately inspired  the next  
development:

\subsection{Construction of Galois representations, part III: Scholze} \label{section:scholze}

In~\cite{scholze-torsion}, Scholze succeeded in constructing Galois representations associated to torsion classes
in the setting of~$\GL_n(\A_F)$ for a CM field~$F$.
  Scholze had the idea after seeing some lectures on~\cite{HLTT}:

\smallskip

\begin{footnotesize}
During a HIM trimester at Bonn, Harris and Lan gave some talks about
their construction of Galois representations. At the time, I had some
ideas in my head that I didn't have any use for: That Shimura varieties
became perfectoid at infinite level, and that there is a Hodge--Tate
period map defined on them. The only consequence I could draw from this
were certain vanishing results for completed cohomology as conjectured
in your work with Emerton; so at least I knew that the methods were able
to say something nontrivial about torsion classes in the cohomology.
After hearing Harris' and Lan's talks, I was trying to see whether these
ideas could help in extending their results to torsion classes. After a
little bit of trying, I found the fake-Hasse invariants, and then it was
clear how the argument would go.

\end{footnotesize}

Even after this breakthrough,  Scholze's construction still fell short of the conjectures in~\cite{CG}
in two ways. The first was that the Galois representation (ignoring here issues of pseudo-representations) was valued not in~$\T$
but in~$\T/I$ for some ideal~$I$ of fixed nilpotence. This is not a crucial obstruction to the methods of~\cite{CG}.
The second issue, however, was that the Galois representations were constructed (in the end) via $p$-adic congruences,
 and thus one did not have
control over their local properties at~$p$ which \emph{are} crucial for modularity applications.

\subsection{The Taylor--Wiles method when  \texorpdfstring{$l_0 > 0$}{l g 0}, part II: DAG}
 
Although not directly related to new~$R=\T$ theorems,  one new recent idea in the subject has been
the work of Galatius--Venkatesh~\cite{GalV} on derived deformation rings in the context of Venkatesh's conjectures over~$\Z$.
This work (in part) reinterprets the arguments of~\cite{CG} in terms of a derived Hecke action.
The authors define a derived version~$\RR$ of~$R$ with~$\pi_0(\RR) = R$.
Under similar hypotheses to~\cite{CG}, the higher
homotopy groups of~$\RR$ are shown to exist precisely in degrees~$0$ to~$l_0$.
One viewpoint of the minimal case of~\cite{CG} is that one constructs a (highly non-canonical) formally smooth ring~$R_{\infty}$
of dimension~$n - l_0$
with an action of a formally smooth ring~$S_{\infty}$ of dimension~$n$ such that the minimal deformation ring~$R$ is~$R_{\infty} \otimes_{S_{\infty}}
S_{\infty}/\a$ for the augmentation ideal~$\a$. Moreover, the ring~$R$ is identified both with the action of~$\T$ on the entire cohomology and simultaneously
on the cohomology in degree~$q_0 + l_0$. 
 On the other hand, when~$l_0 > 1$, the intersection of~$R_{\infty}$ and~$S_{\infty}/\a$ over~$S_{\infty}$ is never
  transverse\footnote{When~$l_0 = 1$, the intersection \emph{can} be transverse when~$R$ is a finite ring. In this setting, the relevant
  cohomology  is also non-zero and finite in exactly one degree.
  On the other hand, as soon as~$\mathrm{Hom}(R,\Qbar_p)$  is non-zero (for example, when there exists an associated motive) and~$l_0 > 0$, the intersection will 
  always be  non-transverse.},
 and homotopy groups of the derived intersection recover the cohomology in all degrees
 (under the running assumptio2n one also knows that the patched cohomology is free).
On the other hand, there is a more canonical way to define~$\RR$, namely to take the unrestricted global deformation ring~$\Rglob$ (which has
no derived structure) and intersect it with a suitable local crystalline deformation ring as algebras over the  unrestricted local deformation ring. The expected dimension
of this intersection is also~$-l_0$ over~$W(k)$, although this is not so clear from this construction. Hence~\cite{GalV}  can be viewed
as giving an intrinsic definition of~$\RR$ independent of any choices of Taylor--Wiles primes and showing that its homotopy groups are related 
(as with~$R_{\infty} \otimes^{\mathbf{L}}_{S_{\infty}} S_{\infty}/\a$) to the cohomology\footnote{There are some subtleties as to what the precise statement
should be in the presence of global congruences, but already this author gets confused at the best of times between homology and cohomology, so I will not try to unentangle these issues here.}. These ideas have hinted at a closer connection between the Langlands program in the arithmetic case and the function field case than was previously
anticipated\footnote{Not anticipated by many people, at least; Michael Harris has been proselytizing the existence of a connection for quite some time.}, see for example work of Zhu~\cite{Zhuarxiv}.

\section{Recent Progress}

\subsection{Avoiding conjectures involving torsion I: the 10-author paper}

As mentioned in~\S\ref{section:scholze}, even after the results of~\cite{scholze-torsion} there remained a significant gap
to make the results of~\cite{CG} unconditional, namely, the conjecture that these Galois representations had the right local
properties at~$p$ and a second conjecture predicting the vanishing of (integral) cohomology localized at a non-maximal
ideal~$\m$ outside a certain precise range (corresponding to known results in characteristic zero). 
It should be noted that the second conjecture was still open (in all but the easiest cases)
in the simpler setting of Shimura varieties. The first hints that  one could possibly make progress on this second conjecture (at least for Shimura
varieties) was given in an informal talk by Scholze in 
Bellairs\footnote{I was invited to give the lecture series in Bellairs after Matthew Emerton didn't respond to his emails. Through
some combination of the appeal of my own work and the fact that the lectures were given on a beach in Barbados,
I managed to persuade Patrick Allen, George Boxer,  Ana Caraiani, Toby Gee, Vincent Pilloni, Peter Scholze, and Jack Thorne to come,
 all of whom are now my co-authors, and all of whom (if they weren't already at the time)
are now more of an expert in this subject than I am. The thought that I managed to teach any of them something about the subject
is pleasing indeed.}  in~2014. This very quickly led to a long term collaboration between Scholze and Caraiani~\cite{CS1,CS2}, which
Caraiani describes as follows:

\smallskip

\begin{footnotesize}
At the Barbados conference in May 2014, Peter gave a lecture on how one might compute the cohomology of compact unitary Shimura varieties with torsion coefficients. The key was to have some control for~$R\pi_{HT*} \mathbb{F}_{\ell}$ restricted to any given Newton stratum. He was expressing this in terms of a conjecture that had grown out of his work on local Langlands using the Langlands--Kottwitz method. After his talk, I went to ask him some questions about this conjecture and it sounded like there were some things that still needed to be made precise. He asked if I wanted to help him make his strategy work. After some hesitation (because I didn’t think I knew enough or was strong enough to work with him), I accepted. Later that evening, I suggested switching from the Langlands--Kottwitz approach to understanding~$R\pi_{HT*}\mathbb{F}_{\ell}$ to an approach more in the style of 
Harris--Taylor. This relies on the beautiful Mantovan product formula that describes Newton strata in terms of Rapoport--Zink spaces and Igusa varieties. Maybe something like this could help illuminate the geometry of the Hodge--Tate period morphism? Peter immediately saw that this should work and we made plans for me to visit Bonn that summer to continue the collaboration.

As Peter and I were finishing writing up the compact case, it became clear to us that the vanishing theorem would give a way to construct Galois representations associated to generic mod~$p$ classes that preserves the desired information at~$p$. Peter started thinking about the non-compact case and how that might apply to the local-global compatibility needed for Calegari--Geraghty. I remember discussing this with him at the Clay Research conference in Oxford in September 2015. By spring 2016, Richard started floating the idea of a working group on Calegari--Geraghty and found out that Peter and I had an approach to local-global compatibility. Around June 2016, Richard suggested to me to organize the working group with him. Peter was very excited about the idea, but wasn’t sure he would be able to attend for family reasons. In the end, we found a date in late October 2016 that worked for everyone.

\end{footnotesize}

The working group met under the auspices of the first ``emerging topics'' workshop\footnote{Although 
later described as
 a ``secret'' workshop, it was an ``invitation-only working group.''}
at the IAS to determine the extent to which the expected consequences could be applied to modularity lifting:
A clear stumbling point
was the vanishing of integral cohomology after localization outside the range of degrees~$[q_0,q_0 + l_0]$. On the other hand,
Khare and Thorne had already observed in~\cite{KT} by a localization argument that this could sometimes be avoided in certain
minimal cases. It was this argument we were able to modify for the general case,
thus avoiding the need to prove the (still open) vanishing conjectures for torsion classes\footnote{I regard my main contribution to~\cite{10author} as explaining how the arguments in~\cite{CG} using Taylor's
Ihara avoidance (\S\ref{section:iharaavoid}) were incompatible with any characteristic zero localization argument in the absence of (unknown)
integral vanishing results in cohomology. The objection (even in the case $l_0 = 0$) was that it was easy to construct complexes~$P^1$
and~$P^2$ of free~$S_{\infty}$ modules so that the support of~$H^*(P^1/p)$ and~$H^*(P^2/p)$ coincided (as they must) but
that (for example) $H^*(P^1)[1/p]$ was zero even though~$H^*(P^2)[1/p]$ was not. The objection to this objection, however,
which was resolved during the workshop (and which to
be clear I played no part in resolving!) is to not merely to compare the support of the complexes~$P^i/p$
but to consider the entire complex in the derived category. In particular, even (say) for a finite~$\Z_p$-module~$M$,
the module~$M[1/p]$ is non-zero exactly when~$M \otimes^{\mathbf{L}} \F_p$ has non-zero Euler characteristic.}.
 The result of the workshop was a success beyond what we could have reasonably 
anticipated
 --- we ended up with
more or less\footnote{It is worth emphasizing that an incredible amount of work was required to turn these ideas into reality,
and that this intellectual effort was by and large
carried out by the younger members of the collaboration.}
 the outline of a plan to prove all the main modularity lifting theorems which finally appeared in~\cite{10author},
namely the Ramanujan conjecture for regular algebraic automorphic forms for~$\GL_2(\A_F)$ of weight zero for any CM field~$F$,
and potential modularity (and the Sato--Tate conjecture) 
 for elliptic curves over CM fields.

There have already been a number of advancements beyond~\cite{10author} including in particular by
Allen, Khare, and Thorne~\cite{AKT} proving the modularity of many elliptic curves over~CM fields and
a potential automorphy theorem for ordinary representations by Qian~\cite{2104.09761}. It does not seem completely implausible
that results of the strength of~\cite{BLGGT} for~$n$-dimensional regular Galois representations of~$G_{\Q}$ 
are within reach.

\subsection{Avoiding conjectures involving torsion II: abelian surfaces}
\label{section:GSp4}
A second example that Geraghty and I had considered during the 2010--2011 IAS special year was the case of abelian surfaces,
corresponding to low (irregular) weight Siegel modular forms of genus~$g = 2$.
It was clear that a key difficulty was proving the vanishing of~$H^2(X,\omega^2)_{\m}$ where~$X$
was a (compactified) Siegel~$3$-fold with good reduction at~$p$, where $\m$ is maximal ideal of the Hecke algebra
corresponding to an absolutely irreducible representation,  and where~$\omega |_{Y} = \det \pi_* \Omega^1_{\AA/Y}$
on the open moduli space~$Y \subset X$ admitting a corresponding
universal abelian surface~$\AA/Y$.
 In other irregular weights (corresponding to motives
with Hodge--Tate weights~$[0,0,k-1,k-1]$ for~$k \ge 4$) the vanishing of the corresponding cohomology groups
was known by Lan and Suh~\cite{LS}. The vanishing of~$H^2(X,\omega^2)_{\m}$ was more subtle, however,
because the corresponding group does \emph{not} vanish in general before localization in contrast to the previous cases.
In~\cite{CGGSp4}, we proved a minimal modularity theorem for these higher weight representations and a minimal
modularity theorem in the abelian case contingent on the vanishing conjecture above which we did not manage to resolve (and which remains unresolved).
I finished and then submitted the paper after I had moved to Chicago and Geraghty had moved to Facebook in 2015.
I then started working with Boxer and Gee\footnote{George Boxer had also arrived at Chicago in 2015, and was collaborating with Gee on companion form
results for Siegel modular forms, with the hope  (in part) of deducing the modularity of abelian surfaces from Serre's
conjecture for~$\mathrm{GSp}_4$ in a manner analogous to the deduction by of the Artin conjecture
from Serre's conjecture for~$\GL_2$ in~\cite{KArtin,KW1}. They usually worked together at Plein Air cafe. Since I had
thought about similar questions with Geraghty and frequently  went to Plein Air for  6oz  cappuccinos, it was not entirely
surprising for us to start working together.} 
on this vanishing question under certain supplementary local hypotheses. 
(By this point,  Galois representations associated to torsion classes in coherent cohomology  had been constructed by Boxer~\cite{1507.05922}
and Goldring--Koskivirta~\cite{Wushi}.)
But then in November of 2016 (one week after the
IAS workshop!), Pilloni's paper on higher Hida theory~\cite{pilloniHidacomplexes} was first posted.
It was apparent to us that Pilloni's ideas would be extremely useful, and the four of us began a collaboration
almost immediately. Just as in~\cite{10author}, we were ultimately able to avoid proving any vanishing conjectures.
However, unlike~\cite{10author}, the way around this problem was not purely by commutative algebra, but instead by
working with ideas from~\cite{pilloniHidacomplexes}. Namely, instead of working with the  cohomology of the full Siegel modular variety~$X$,
one could work with
 the coherent cohomology of a certain open variety of~$X$ with cohomological dimension one whose (infinite dimensional)
cohomology could still be tamed using the methods of higher Hida theory~\cite{pilloniHidacomplexes}
in a way analogous to how Hida theory controls the (infinite dimensional) cohomology of the affine variety (with cohomological dimension zero)
corresponding to the ordinary locus. Generalizing this to a totally real field, one could then \emph{combine} these
ideas with the Taylor--Wiles method as modified in~\cite{CG}
to prove the potential modularity of abelian surfaces over totally real fields~\cite{BCGP}.
This coincidentally gives a second proof of the potential modularity of elliptic curves over CM fields proven in~\cite{10author}.
(The papers~\cite{10author} and~\cite{BCGP}    both were conceived of and completed
within a week or so of each other.)

\section{The depths of our ignorance}
\label{section:hopeless}
Despite what can reasonably be considered significant progress in proving many cases of modularity since 1993, it remains the case that
many problems appear just as hopeless as they did then\footnote{Or perhaps harder, since there has been almost 30 years without
any progress whatsoever.}. Perhaps most embarrassing is the case of even Galois representations~$G_{\Q} \rightarrow \GL_2(\C)$ with
non-solvable image (equivalently, projective image~$A_5$). For example, we cannot establish the Artin conjecture for a \emph{single} 
Galois representation whose image is the binary icosahedral group~$\SL_2(\F_5)$ of order~$120$. 
The key problem is that the automorphic forms (Maass forms with
eigenvalue~$\lambda = 1/4$ in this case) are very hard to access --- given an even (projective)~$A_5$ Galois
representation, we don't even know how to prove that there exists a corresponding Maass form with the right Laplacian eigenvalue,
let alone one whose Hecke eigenvalues correspond to the Galois representation\footnote{Motives can be divided according to a tetrachotomy.
The first form are the Tate (and potentially Tate) motives, whose automorphy was known to Riemann and Hecke.
The second form are the motives (conjecturally) associated to automorphic representations which
are discrete series at infinity and thus amenable to the Taylor--Wiles method.
The third form are the motives (conjecturally) associated to automorphic representations
which are at least seen by some flavour of cohomology, either by  the Betti cohomology of locally symmetric spaces or the coherent cohomology
of Shimura varieties (possibly in degrees greater than zero) which are amenable in principle to the modified Taylor--Wiles method. The fourth form consist
of the rest, which (besides a few that can be accessed
by cyclic base change) are a complete mystery.}.
 In many ways, we have made no real progress on this question.
The case of curves of genus~$g > 2$ whose Jacobians have no extra endomorphisms seems equally hopeless.
One can only take solace in the fact that the Shimura--Taniyama conjecture seemed equally out of reach before Wiles' announcement in Cambridge
in 1993.

\bibliographystyle{emss}
\bibliography{fcalegari}

\section{Acknowledgments}

The arithmetic Langlands program has been a bustling field of activity in the past~$30$ years with many deep
results by a number of extraordinary mathematicians. To be part of this field (mostly as an observer but
occasionally as a contributor) has been an amazing experience.
I would  like to thank my collaborators all of whom have taught me so much, with special thanks to Matthew Emerton
and Toby Gee who have in addition been my close friends.
Thanks to Toby Gee and Ravi Ramakrishna for a number of useful suggestions and corrections on the first version of this survey.
Finally, 
special thanks to those mathematicians who allowed me to share some of their personal recollections in this paper,
in particular Christophe Breuil,  Ana Caraiani, Pierre Colmez,
Fred Diamond, Michael Harris, Mark Kisin, 
Peter Scholze, Richard Taylor,
and Andrew Wiles.

\end{document}